\documentclass[journal,oneside,web,onecolumn]{ieeecolor}
\usepackage{lcsys}

\usepackage{graphicx} 
\usepackage{graphics}
\usepackage{textcomp}
\usepackage{epsfig}
\usepackage{times} 
\usepackage{amsmath}
\usepackage{amssymb}
\usepackage{amsfonts}
\usepackage{float}
\usepackage{siunitx}
\usepackage{scalerel}
\usepackage{cite}
\def\BibTeX{{\rm B\kern-.05em{\sc i\kern-.025em b}\kern-.08em
    T\kern-.1667em\lower.7ex\hbox{E}\kern-.125emX}}

\newtheorem{theorem}{Theorem}
\newtheorem{proposition}{Proposition}
\newtheorem{lemma}{Lemma}

\newtheorem{remark}{Remark}
 \newtheorem{assumption}{Assumption}
 \newtheorem{Corollary}{Corollary}

\title{\LARGE \bf Distributed Estimation}
\newcommand{\ncom}{\newcommand}

\newcommand{\beqn}{\begin{eqnarray*}}
\newcommand{\eeqn}{\end{eqnarray*}}
\newcommand{\beq}{\begin{eqnarray}}
\newcommand{\eeq}{\end{eqnarray}}

\newcommand{\norm}[1]{\left\lVert #1 \right\rVert}
\newcommand{\inprod}[2]{\left\langle #1, #2 \right\rangle}
\ncom\R{\mathbb{R}}
\DeclareMathOperator*{\argmin}{arg\,min}
\DeclareMathOperator*{\argmax}{arg\,max}

\DeclareMathOperator*{\bd}{bd}
\DeclareMathOperator*{\Lim}{Lim}

\DeclareMathOperator*{\Liminf}{Liminf}
\DeclareMathOperator*{\dom}{dom}

\author{Anik Kumar Paul, Arun D Mahindrakar and Rachel K Kalaimani
\thanks{Anik  is a Graduate student in the Department of Electrical Engineering, IIT Madras, Chennai-600036, India
        {email: anikpaul42@gmail.com}}
        \thanks{Arun  and Rachel are  with the Department of Electrical Engineering, Indian Institute of Technology Madras, Chennai-600036, India (email: arun\_dm@iitm.ac.in, rachel@ee.iitm.ac.in) }
    }

\begin{document}
\title{Convergence Analysis of Stochastic Saddle Point Mirror Descent Algorithm - A Projected Dynamical View Point}
\maketitle
\thispagestyle{plain}
\pagestyle{plain}
\begin{abstract}
Saddle point problems, ubiquitous in optimization, extend beyond game theory to diverse domains like power networks and reinforcement learning.  This paper presents novel approaches to tackle saddle point problem, with a focus on continuous-time contexts. In this paper we propose a continuous time dynamics to tackle saddle point problem utilizing projected dynamical system in non-Euclidean domain. This involves computing the (sub/super) gradient of the min-max function within a Riemannian metric. Additionally, we establish viable Caratheodory solutions also prove the  Lyapunov stability and asymptotic set stability of the proposed continuous time dynamical system.
Next, we present the Stochastic Saddle Point Mirror Descent (SSPMD) algorithm and establish its equivalence with the proposed continuous-time dynamics. Leveraging stability analysis of the continuous-time dynamics, we demonstrate the almost sure convergence of the algorithm's iterates. Furthermore, we introduce the Zeroth-Order Saddle Point Mirror Descent (SZSPMD) algorithm, which approximates gradients using Nesterov's Gaussian Approximation, showcasing convergence to a neighborhood around saddle points. The analysis in this paper provides geometric insights into the mirror descent algorithm and demonstrates how these insights offer theoretical foundations for various practical applications of the mirror descent algorithm in diverse scenarios.
\end{abstract}


\section{Introduction}
Saddle point problem, popularly termed as a min-max problem, entails an optimization problem of an objective function characterized by two parameters. In this scenario, the objective  is to minimize the given  objective function with respect to one variable while concurrently maximizing it with respect to another variable.
 Saddle point optimization predominantly involves identifying Nash equilibrium within zero-sum games, and has long been a a central focus of extensive research within the field of game theory \cite{9833083,9833242}.  However, the scope of saddle point problems and the algorithms used to tackle them go beyond game theory, finding applications in a myriad of domains, including but not limited to power networks, adversarial training, reinforcement learning, and adversarial networks \cite{NIPS2005e465ae46, pmlr-v80-dai18c,  du2017stochastic, NIPS200464036755}.   The exhaustive analysis and performance analysis of diverse algorithms for tackling the saddle point problem have constituted a substantial area of scholarly interest over an extended period, encompassing both discrete time  \cite{saddle,Zhao2019OptimalAF} and continuous time \cite{8123914,7171030,doi:10.1137/15M1026924}. In practical applications, saddle point algorithms are commonly implemented in discrete time; however, their continuous-time analysis remains a focal point of considerable interest within the realms of control and optimization. The implementation of the saddle point algorithm often involves discretizing its continuous-time counterpart, and the convergence analysis of the algorithm can be effectively done through the application of the theory of differential equations.  For example, in \cite{doi:10.1137/15M1026924}, the analysis focuses on stability by studying a set of ordinary differential equations (ODEs) known as saddle point dynamics. These ODEs capture the gradients of the max-min function with respect to both the variables. Investigating the stability of these equations provides a solid foundation for understanding how the saddle point algorithm behaves and converges over time in continuous time context. However, the analysis in \cite{doi:10.1137/15M1026924} centers on showing the stability of saddle point dynamics, under the assumption that the saddle point problem is unconstrained. In the context of a constrained saddle problem, the analysis necessitates the consideration of a projected dynamical system. This framework ensures under mild assumptions that the solutions derived from the ODEs remain confined within the desired sets (constrained set) imposed by the saddle problem. The analysis of a projected dynamic system introduces a degree of complexity beyond that of standard ODEs. This intricacy arises due to the presence of a discontinuous vector field in the right-hand side (RHS) of the projected dynamical system. In  \cite{8123914,7171030},  projected dynamical system is employed to address constrained optimization problems featuring both equality and inequality constraints.  The constrained optimization problem can be represented equivalently as saddle point problems, with a constraint on the dual variable in view  of  Karush-Kuhn-Tucker (KKT) conditions (see \cite{citeulike:163662} for more details). For a broader perspective on the saddle point problem, the existence and stability of projected dynamical systems have been discussed in \cite{GOEBEL201716}.

As highlighted in \cite{doi:10.1137/18M1229225}, the existing projected dynamical system in the Euclidean domain, though effective for modeling many constrained optimization problems, has limited applicability when dealing with complex optimization scenarios. Notably, it overlooks the geometric aspects of constrained optimization problems. Motivated by the insights from the analysis presented in \cite{doi:10.1137/18M1229225}, our first contribution in this paper is the development of a projected dynamical system in a non-Euclidean domain to address the saddle point problem.  We formulate the (sub/super) gradient of the min-max function within a Riemannian metric framework and compute the projections with respect to this Riemannian metric.  This formulation significantly broadens the scope for tackling a wider array of complex problems.  It's important to note that the min-max function considered here is nonsmooth convex-concave, leading to a more general projected differential inclusion to tackle the saddle point problem.
Consequently, the subsequent contribution of this paper is to extend the existence result of Caratheodory solution provided in \cite{doi:10.1137/18M1229225} for the projected differential inclusion. 
Further, we demonstrate Lyapunov stability of all solutions of the proposed dynamics with respect to the set of saddle points as equilibrium points, along with asymptotic set stability to the set of saddle points.

However, an important question arises: despite the enhanced geometric understanding provided by the proposed projected dynamical system, how is its theoretical analysis beneficial in practical real-life situations? Additionally, how can we approximate it in the discrete time domain for practical implementation? For instance, in \cite{projected}, the approximation of the projected dynamical system in Euclidean domain for numerical methods is elucidated, resulting in the formulation of the projected subgradient method for convex optimization problems.
The next contribution of this paper is to demonstrate that the proposed projected dynamical system in the Riemannian manifold can be numerically approximated as iterate steps of the standard mirror descent algorithm.   Originally, the mirror descent algorithm generalizes the standard gradient descent algorithm defined over Euclidean space and applies to a optimization problem over Hilbert and Banach spaces \cite{8409957}. This extension establishes a versatile framework, facilitating optimization in spaces characterized by intricate geometries or the inclusion of additional structural complexities.  For example, it is known that for many constrained sets like probability simplex or, $p$- norm ball, the implementation of mirror descent algorithm is much more computationally efficient than the standard projected subgradient method.
Recent years have witnessed a surge in the prominence of the mirror descent algorithm, particularly in the realms of large-scale optimization problems, data-driven control and learning, power systems, robotics, and game theoretic applications \cite{8409957,anik}. From the discussion in the paper, it becomes evident that the projected dynamical system in the Riemannian manifold can be numerically implemented through the steps of the saddle point mirror descent algorithm.

In the next section of the paper, we explore SSPMD algorithm. In various practical scenarios, including robust estimation \cite{1055386}, resource allocation, and learning theory \cite{KellyMT98,811451}, the objective function is influenced by inherent randomness. In such cases, computing exact gradients for the min-max function can prove impractical or computationally intensive. This provides the motivation to consider the SSPMD algorithm.
The SSPMD algorithm was initially introduced in \cite{rsa}, where the convergence analysis of the function value of the iterates of the SSPMD algorithm was shown in expectation. However, in this paper, we focus on the almost sure convergence of iterates to the set of saddle points. An almost sure convergence guarantee  is more significant than the convergence in expectation since it describes what happens to the individual trajectory in each iteration. The proof strategy for establishing almost sure convergence of the iterates of the SSPMD algorithm relies on the dynamic viewpoint of stochastic approximation, as discussed in \cite{doi:10.1137/S0363012904439301}. The core idea of this dynamic viewpoint is that the continuous-time interpolation of any iterates stochastic recursive algorithm can be seen as a ``pseudo-trajectory" of a limiting ordinary differential equation or inclusion. This concept has been extensively used to demonstrate almost sure convergence analysis of projected gradient descent, where the limiting differential equation is proven to be the projected dynamical system in the Euclidean domain (refer to \cite{borkar2008stochastic} and \cite{kushner1978stochastic} for more details).  This concept is further extended to demonstrate almost sure convergence analysis of Nesterov's dual averaging or weakly mirror descent by showcasing the stability analysis of a suitable limiting differential equation. The primary question in demonstrating the almost sure convergence of the stochastic mirror descent algorithm is determining the corresponding limiting differential equations or inclusions. In this paper, we show that this limiting differential inclusion is the projected differential inclusion in the non-Euclidean domain mentioned earlier. By leveraging the stability analysis of projected differential inclusions, we propose a novel proof strategy to demonstrate the almost sure convergence of the SSPMD algorithm.  It is worth noting that in  \cite{doi:10.1137/120894464} , the almost sure convergence analysis of the stochastic mirror descent algorithm is demonstrated by employing a probabilistic approach, primarily centered around the theory of martingales and the associated concept of almost supermartingales. This probabilistic approach can be extended to demonstrate the almost sure convergence of the SZSPMD  algorithm. While the choice between these approaches ultimately depends on personal preference, focusing on the differential inclusion method provides a comprehensive analytical framework for the algorithm.
Here, we emphasize the differential inclusion approach as it introduces a dynamic perspective that offers deeper insights into the algorithm compared to the standard analysis of its discrete-time counterpart using martingale theory. Consequently, the convergence analysis of the algorithm can be carried out using nonlinear system theory.  This approach proves to be particularly helpful in many applications where there are errors in computing the gradient of the function. This is further demonstrated in the next contribution, where we present the SZSPMD algorithm.

For the SZSPMD algorithm, we presume the existence of an oracle that can provide the function value at any given point without disclosing information about the gradient (assuming the min-max function is continuously differentiable). Each step in this algorithm resembles that of the stochastic saddle point mirror descent, with the key difference being the need to approximate the function's gradient at each point. This approximation relies on Nesterov's Gaussian Approximation which depends  on a smoothness parameter, as outlined in \cite{Nesterov2015RandomGM}.  There has recently been a surge of interest generated in different variants of  zeroth-order optimization, for both convex and non-convex functions (see \cite{9186148} and the references therein).  For the SZSPMD algorithm, we have showcased that the limiting differential inclusion takes the form of a projected differential inclusion with a disturbance input.  In other words, the continuous-time interpolation of the iterates of the SZSPMD   algorithm emerges as an asymptotic pseudo trajectory of the projected dynamical system in a Riemannian manifold with a disturbance input.  Through robust analysis, we show that by regulating the magnitude of the disturbance input (controllable via the choice of smoothing parameter), the iterates of the zeroth-order saddle point mirror descent algorithm almost surely converge to any desired neighborhood of the set of saddle points. This subsection's contribution extends the findings presented in \cite{8016343}. In that work, they established the almost sure convergence of stochastic gradient descent using the perspective of differential inclusions. However, that paper specifically dealt with unconstrained optimization problems.
In this paper, we do not make any such assumptions. Since, the analysis in this paper extends beyond the traditional Euclidean domain, the outcomes detailed in this paper also apply to the zeroth-order stochastic gradient descent algorithm in the context of constrained optimization problems.

The following list summarises the key contribution of this
study.

\begin{enumerate}
    \item We present a continuous-time dynamics approach using a projected dynamical system in a non-Euclidean space to tackle the  non-smooth saddle point problem. This method involves computing the (sub/super) gradient of the min-max function within a Riemannian metric framework and projecting based on this metric. Due to the inherent non-smoothness of the saddle point problem, the right-hand side of the continuous-time dynamics results in set-valued maps, giving rise to a projected differential inclusion in a non-Euclidean domain, which contrasts with the usual projected differential equation in non-Euclidean domain found in most recent  literatures. 
   We establish the existence of a viable Caratheodory solution for this projected dynamical system, represented equivalently through a normal cone representation. Additionally, we show that the saddle points coincide with the equilibrium points of this dynamical system. Furthermore, we prove that all saddle points serve as Lyapunov stable equilibrium points, and all Caratheodory solutions of the projected dynamical system asymptotically converge to the set of  saddle points. 
   
\item We propose numerical approximation techniques for the projected dynamical system in the non-Euclidean domain, facilitating practical implementation through iterative steps of the mirror descent algorithm. In this paper, the convergence analysis is conducted on a more generalized version - the SSPMD algorithm, where exact (sub/super) gradient computation of the min-max function is computationally challenging. The emphasis in this paper lies in showing almost sure convergence of the iterates of the SSPMD algorithm. Leveraging the stability analysis of projected differential inclusion, we establish that the iterates of the SSPMD algorithm almost surely converge to the set of saddle points.
\item We finally introduce  SZSPMD algorithm where the gradient of the min-max function is approximated using Nesterov's Gaussian Approximation. Here, we demonstrate that the continuous time interpolation of the iterates of the SZSPMD algorithm can be approximated with the Caratheodory solution of projected dynamical system in non-Euclidean domain, albeit with a disturbance input. Through robust analysis, we demonstrate that iterates of the SZSPMD algorithm converge to a neighborhood around the set of saddle points, with the size of this neighborhood being influenced by smoothing parameters of the approximated gradient.
\end{enumerate}

\section{Notation and Math Preliminary}
\label{mathpre}
\subsection{Notation}
Let $\R$, $\mathbb{R}^{+}$ and $\mathbb{R}^n$ represent the set of real numbers, set of non-negative numbers,  set of $n$ dimensional real vectors.   Let~$\norm{.}$ denote any \mbox{norm} on~$\R^n$.  Given a norm $\norm{.}$ on $\mathbb{R}^n$, the dual norm of $x\in \R^n$ is $\norm{x}_\ast := \sup\{ \inprod{x}{y}: {\norm{y}\leq 1}, y\in \R^n \}$. The notation $x[i]$ denotes the $i$th component of the vector $x$. $I_n$ is $n \times n$ identity matrix. 
If $x \in \mathbb{R}^n$, then $U(x,\epsilon)$ represents  an open ball in $\mathbb{R}^n$ centred at $x$ with radius $\epsilon$.  For any set $\mathcal{X}$,  $\bd (\mathcal{X})$ denotes its boundary. Let $\mathcal{X}$ be a closed set in $\mathbb{R}^n$ then $\mathcal{I}_\mathcal{X}(x)$ denotes the indicator function of the set $\mathcal{X}$.
Consider $L(x,\gamma) : \mathbb{R}^n \times \mathbb{R}^m \to \R$ is continuously differentiable at each $(x,\gamma) \in \mathbb{R}^n \times \mathbb{R}^m \to \mathbb{R} $ then $\nabla_x L(x,\gamma) = (\frac{\partial L}{\partial x_1} , \frac{\partial  L}{\partial x_2}  , \ldots, \frac{\partial L}{\partial x_n})$. Similarly, define $\nabla_\gamma L(x,\gamma)$.
A random vector  $X\sim \mathcal{N}(0_n, I_n)$ denotes a $n$-dimensional normal random vector with zero-mean and unit standard-deviation. For two random variables $X$ and $Y$,  $\sigma(X,Y)$ is the smallest sigma-algebra generated by  random variables $X$ and $Y$. Let $(\{X_i\}, \; 1 \leq i \leq n )$ be the sequence of random vectors defined on the probability space $(\Omega,\mathcal{F},\mathbb{P})$ then, $\sigma(\{X_i\}, \; 1 \leq i \leq n )$ denotes the smallest sigma algebra generated by the random vectors $(\{X_i\}, \; 1 \leq i \leq n)$.
\subsection{Math Preliminaries}
Let $\mathcal{X}\subseteq \R^n$ be a convex set and $f: \mathbb{R}^n \to \Bar{\mathbb{R}}$  a convex function. Throughout our analysis, we assume that $f$ is proper, meaning its effective domain is nonempty, denoted by
\begin{equation*}
    \dom{f} (\subseteq \mathcal{X}) = \{x \in \mathbb{R}^n \; | \; f(x) < \infty \} \neq \emptyset.
\end{equation*}

A vector $g\in \R^n$ is called a subgradient of $f$ at $x$ if
    \begin{equation*}
        f(y)\ge f(x)+\inprod{g}{y-x}\;\ \; \forall \; y\in \R^n .
    \end{equation*}
     The set of all subgradients of $f$ at $x$ is termed as subdifferential of $f$ at $x$, denoted by $\partial f(x)$. If $f$ is  convex and continuously differentiable at $x\in \R^n$, then $\partial f(x)=\{\nabla f(x)\}$, where $\nabla f(x)$ is the gradient of $f$ at $x$.
It is well known property that the subdifferential  of convex function  is nonempty, convex and compact.  If $f$ is $L$-Lipschitz over a compact set $\mathcal{X}$, then $\sup\limits_{g \in \partial f (x)} \norm{g}  \leq L$. Consider any $x \in \mathcal{X}$, then the directional derivative of $f$ at $x$ in the direction of $v \in \mathbb{R}^n$, denoted by $D f(x)[v]$, is defined as 
\begin{equation*}
    D f(x) [v] = \lim\limits_{t \downarrow 0} \frac{f(x+tv)-f(x)}{t}.
\end{equation*}

The subdifferential $\partial f(x)$ can also be characterized by the following equation

\begin{equation*}
    \partial f(x) = \{ g \in \mathbb{R}^n \; | \inprod{g}{v} \leq  D f(x) [v] \; \; \forall \; v \in \mathbb{R}^n \}. 
\end{equation*}

Let  $R : \mathbb{R}^n \to \Bar{\mathbb{R}}$ be differentiable over an open set that contains $\mathcal{X}$. The function $R$ is said to be $\sigma_R$-strongly convex over $\mathcal{X}$  if, for all $x,y \in \mathcal{X}$, 
  \begin{equation*}
      R(y) \geq R(x) + \inprod{\nabla{R}(x)}{y-x}+ \frac{\sigma_R}{2}\norm{y-x}^2.
  \end{equation*}
  Furthermore, for any $x,y \in \mathcal{X}$, the following inequality holds  \cite{8187324}
  \begin{equation*}
      \inprod{\nabla R (x) - \nabla R(y)}{x-y} \geq \sigma_R \norm{x-y}^2 .
  \end{equation*}
The function  $f$ is termed lower semicontinuous at any $x \in \mathcal{X}$ if for any sequence $\{x_n\} \subseteq \mathcal{X}$ satisfying $\lim\limits_{n \to \infty} x_n = x$, it holds that  
$f(x) \leq \liminf\limits_{n \to \infty} f(x_n)$. Since the function $f$ is convex, we can equivalently say $\dom f$ is closed in $\mathbb{R}^n$.

The function $h: \mathbb{R}^n \to \mathbb{R} \cup \{-\infty\}$ is concave and uppersemicontinuous if $-h$ is convex and lower semicontinuous.

 The conjugate of the function $f$, denoted by $f^\ast$, is defined as
\begin{equation}
    f^\ast(z) = \sup\limits_{x \in \mathbb{R}^n} \{ \inprod{x}{z}-f(x) \}.
\end{equation}
In addition, if $f$ is strongly convex then $f^\ast: \mathbb{R}^n \to \mathbb{R}$ is continuously differentiable, where
\begin{equation*}
    \nabla f^\ast(z) = \argmax\limits_{x \in \mathbb{R}^n} \{ \inprod{x}{z}-f(x) \}.
\end{equation*}
    Suppose $f$ is differentiable over an open set that contains $\mathcal{X}$. Then, according to first-order optimality condition, a vector $x^\ast \in \mathcal{X}$ minimizes the function $f$ over  $\mathcal{X}$  if and only if
  \begin{equation*}
      \nabla f(x^\ast)^\top (x-x^\ast) \geq 0 \;\;  \forall \; x \in \mathcal{X}.
  \end{equation*}
Let $\mathcal{A}$ be a closed connected subset of $\mathbb{R}^n$, and let $x$ be any point in $ \mathbb{R}^n$. The distance between $x$ and $\mathcal{A}$ is defined as
\begin{equation*}
    d(x,\mathcal{A}) = \inf\limits_{a \in \mathcal{A}} \norm{a - x}.
\end{equation*}
The $\epsilon$- neighborhood of the set $\mathcal{A}$, denoted  $N_\epsilon(\mathcal{A})$, is defined as
\begin{equation*}
    N_\epsilon (\mathcal{A}) = \{ x \in \mathbb{R}^n \; | \; d(x,\mathcal{A}) \leq \epsilon \}.
\end{equation*}

 Consider a closed convex set $\mathcal{X} \subseteq \mathbb{R}^n$. The tangent cone to the set $\mathcal{X}$ at a point $x \in \mathcal{X}$ is defined as follows
 \begin{equation*}
     \mathcal{T}_\mathcal{X}(x) = \{ v \in \mathbb{R}^n \; | v= \lim\limits_{n \to \infty} \frac{x_n-x}{t_n}, \; \;   x_n \in \mathcal{X}\; \; \text{and} \;  t_n \downarrow 0 \}. 
 \end{equation*}
 
 It can be verified that $\mathcal{T}_\mathcal{X}(x) = \Liminf\limits_{x_n \to x} \mathcal{T}_\mathcal{X}(x_n)$. In other words, $\nu \in \mathcal{T}_\mathcal{X}(x)$ if and only if for any sequence $x_n (\in \mathcal{X})$ converges to $x $, $\exists \; \nu_n \in \mathcal{T}_\mathcal{X}(x_n)$ such that $\nu_n \to \nu$.  
 
 Next,  the normal cone to the set $\mathcal{X}$ at $x \in \mathcal{X}$ is the polar of the tangent cone at $x$
 \begin{equation*}
     \mathcal{N}_\mathcal{X}(x) = \mathcal{T}_\mathcal{X}(x)' = \{\eta \in \mathbb{R}^n \; | \inprod{\eta}{v} \leq 0 \; \; \forall \; v \in \mathcal{T}_\mathcal{X}(x) \}.
 \end{equation*}
In this context, let $F: \mathcal{X} \rightrightarrows \mathbb{R}^n$ denote a set-valued map. We consider the following ordinary differential inclusions:
 \begin{equation}
 \begin{split}
     & \Dot{x}(t) \in F(x(t))
     \\ \; \; & x(0) = x_0.
     \end{split}
     \label{odo}
 \end{equation}
 We focus on viable Caratheodory solutions of the ordinary differential inclusion\eqref{odo}. Suppose $x(t) : \mathbb{R}^{+} \to \mathcal{X}$ is an absolutely continuous function satisfying $\dot{x}(t) \in F(x(t))$ almost everywhere with the initial condition $x(0) = x_0$, then $x(t)$ is referred to as a viable Caratheodory solution of the differential inclusion \eqref{odo}.

The following proposition  \cite{aubin1993differential} provides conditions for the existence of viable Caratheodory solutions for the ordinary differential inclusions \eqref{odo} for all initial conditions $x_0 \in \mathcal{X}$.
\begin{proposition}{\cite{4518905}}
    Consider the following conditions:  
    \begin{enumerate}
        \item For all $x \in \mathcal{X}$, $F(x)$ is  compact and convex. Additionally, $F(x)$ is upper semicontinuous at each $x \in \mathcal{X}$, meaning that  $\forall \; \epsilon >0$, $\exists \; \delta > 0$ such that $\forall \; y \in U(x,\delta)$, $F(y) \subseteq F(x) + U(0,\epsilon)$.
        \item For all $x \in \mathcal{X}$,
\begin{equation*}
\sup_{y \in F(x)} \norm{y} \leq \mathrm{K}(1+ \norm{x}),
\end{equation*}
where $\mathrm{K} > 0$.
       \item For all $x \in \mathcal{X}$, $F(x) \cap \mathcal{T}_\mathcal{X}(x) \neq \emptyset$.
    \end{enumerate}
    Then, there exists a viable Caratheodory solution of the ordinary differential inclusion \eqref{odo} for all initial conditions $x_0 \in \mathcal{X}$. 
    \label{dh}
\end{proposition}

Any vector $x_e \in \mathcal{X}$ is an equilibrium point if and only if $0 \in F(x_e)$.

 A set $M \subseteq \mathcal{X}$ is considered weakly invariant for \eqref{odo} if, for every $x_0 \in M$, there exists at least one Caratheodory solution $x(t)$ of \eqref{odo} with $x(0) = x_0$ satisfying $x(t) \in M$ for all $t \geq 0$ almost everywhere.

The set $M$ is called strongly invariant if for every $x_0 \in M$ and  all Caratheodory solutions $x(t)$ of \eqref{odo} with initial condition $x(0) = x_0$ satisfy $x(t) \in M$ for all $t \geq 0$ almost everywhere. 

 Suppose  $V : \mathbb{R}^n \to \mathbb{R}$ and a set valued map $F : \mathcal{X} \rightrightarrows \mathbb{R}^n$, the set valued Lie derivative of $V$ with respect to $F$ at $x \in \mathcal{X}$ is defined as
\begin{equation*}
    \Tilde{\mathcal{L}}_F V(x) = \{ a \in \mathbb{R}\; | \exists \; \; v \in F(x) \; \; \text{such that} \; \; \zeta^\top v = a \; \; \forall \; \zeta \in \partial V(x) \}.
    \end{equation*}
    If $V$ is continuously differentiable at $x \in \mathcal{X}$ then 
    \begin{equation*}
      \Tilde{\mathcal{L}}_F V(x)  = \{ \nabla V(x)^\top v \; | v \in F(x) \}  .
    \end{equation*}
    
  An equilibrium point $x_e \in \mathcal{X}$ of \eqref{odo} is  Lyapunov stable if $\forall \; \epsilon > 0$, $\exists \; \delta(\epsilon) > 0$ such that $\forall \; x_0 \in U(x_e,\delta)\cap \mathcal{X}$, and for all viable Caratheodory solutions of \eqref{odo} $x(t)$ with initial condition $x(0) = x_0$, it holds $x(t) \in U(x_e,\epsilon) \cap \mathcal{X}$ $\forall \; t \geq 0$. 

The following theorem provides sufficient conditions for the stability of \eqref{odo}.
  \begin{theorem}
       Let $x_e$ be an equilibrium point of the differential inclusion \eqref{odo}.  Let $\mathcal{D} \subseteq \mathbb{R}^n$  is an open  and connected set  that contains the point $x_e$.  Furthermore, assume that $V : \mathcal{X} \to \mathbb{R}$ is a continuously differentiable function with $V(x_e) = 0$, and $V(x) >0 \; \forall \; x \in \mathcal{X} \cap \mathcal{D} \setminus
       \{x_e\}$. If $\max \Tilde{\mathcal{L}}_F V(x)  \leq 0$ for each $x \in \mathcal{X} \cap \mathcal{D}$, then $x_e$ is a Lyapunov stable equilibrium point.
       \\ Additionally, if $M \subseteq \mathbb{R}^n$  is a  compact set  and $M \cap \mathcal{X}$ is a  strongly invariant set for the differential inclusion \eqref{odo} such that $\max \Tilde{\mathcal{L}}_F V(y)  \leq 0$ $\forall \; y \in M \cap \mathcal{X}$, then all Caratheodory solutions of the differential inclusions \eqref{odo} starting from $x_0 \in M \cap \mathcal{X}$ converge to the largest weakly invariant set  contained in $\mathcal{X}  \cap M \cap  \{ y \in \mathcal{X} | 0 \in \Tilde{\mathcal{L}}_F V(y) \} $.
       \label{Cortes}
   \end{theorem}
   
   \textbf{Stochastic Approximation:} Let ${x(n)}$ represent a discrete-time stochastic process in $\mathbb{R}^n$, governed by the following recursive equation:

\begin{equation*}
x(n+1) = x(n) + \alpha(n) (y(n) + M(n+1)) + b(n).
\end{equation*}

We assume the following conditions hold:
\begin{itemize}
    \item $y(n) \in F(x(n))$, where $F: \mathcal{X} \rightrightarrows \mathbb{R}^n$ satisfies all conditions specified in Proposition \ref{dh}.
    \item The step-size $\alpha(n)$ satisfies  $\sum\limits_{n \geq 1} \alpha(n) = \infty$ and $\sum\limits_{n \geq 1} \alpha(n)^2 < \infty$.
    \item ${M(n+1)}$ is a martingale difference sequence. This means that for ${\mathcal{F}_n}$, the filtration generated by the random vectors ${x(n)}$ ($\mathcal{F}_n = \sigma ({x(i)}, 1 \leq i \leq n)$), we have $\mathbb{E}[M(n+1) | \mathcal{F}_n] = 0$, and  also $\mathbb{E}[\norm{M(n+1)}^2| \mathcal{F}_n] \leq \mathcal{K}(1 + \norm{x(n)}^2)$, where $\mathcal{K} > 0$.
    \item  $b(n) = \small{o}(\alpha(n))$, meaning that $\lim\limits_{n \to \infty} \frac{\norm{b(n)}}{\alpha(n)} = 0$.
\end{itemize}
Next, we define continuous-time interpolated trajectories $\Hat{x}(t)$ corresponding to the iterates $x(n)$ as follows:
\begin{equation*}
\Hat{x}(t) = x(n) + (x(n+1)-x(n)) \frac{t-t(n)}{t(n+1)-t(n)}, \quad \forall \; t \in I_n,
\end{equation*}
where $I_n = [t(n),t(n+1)]$. According to \cite{doi:10.1137/S0363012904439301}, $\hat{x}(t)$ is an asymptotic pseudotrajectory of the ordinary differential inclusion $\dot{x}(t) \in F(x)$. In other words , for any $T > 0$:

\begin{equation*}
\lim\limits_{t \to \infty} \inf\limits_{x(s) \in S_{\hat{x}(t)}} \sup\limits_{0 \leq s \leq T} \norm{\Hat{x}(t+s)-x(s)} = 0,
\end{equation*}

where $S_{\Hat{x}(t)}$ is the set of all Caratheodory solutions $x(s)$ of the ordinary differential inclusion $\dot{x}(t) \in F(x)$ with the initial condition $x(0) = \Hat{x}(t)$.
\section{Problem Statement}
This paper delves into the analysis of a convex-concave saddle point problem represented as follows:
\begin{equation}
    \min\limits_{x \in \mathcal{X}} \max\limits_{\gamma \in \mathcal{Y}} L(x,\gamma).
    \label{primal}
\end{equation}

Here, the constraint sets $\mathcal{X} \subseteq \mathbb{R}^n$ and $\mathcal{Y} \subseteq \mathbb{R}^m$ are non-empty convex and compact. The function $L: \mathbb{R}^n \times \mathbb{R}^m \rightarrow \overline{\mathbb{R}} \cup \{- \infty\}$ is convex with respect to  the variable $x \in \mathbb{R}^n$ and concave with respect to  the variable $\gamma \in \mathbb{R}^m$. The solutions to the problem \eqref{primal}  are known as saddle points.
These are characterized by the following inequality:

\begin{equation*}
    L(x^\ast,\gamma) \leq L(x^\ast, \gamma^\ast) \leq L(x, \gamma^\ast) \; \; \forall \; x \in \mathcal{X} \; \; \text{and} \; \;  \forall \; \gamma \in \mathcal{Y}.
\end{equation*}

The set of saddle points, denoted by $\mathcal{S}$, is defined as:

\begin{equation*}
    \mathcal{S} = \{(x^\ast, \gamma^\ast) \in \mathcal{X} \times \mathcal{Y} | \; \;  (x^\ast,\gamma^\ast) \text{ is a saddle point} \}.
\end{equation*}

Before proceeding further, we make the following assumptions.

\begin{assumption}
    The set $\mathcal{S}$ is both nonempty and compact in $\mathbb{R}^n \times \mathbb{R}^m$. Additionally,  $\forall \; (x,\gamma) \in (\mathcal{X} \times \mathcal{Y}) \setminus \mathcal{S}$, $\exists \; (x^\ast, \gamma^\ast) \in \mathcal{S}$ such that either $L(x^\ast,\gamma) < L(x^\ast,\gamma^\ast)$ or, $L(x^\ast,\gamma^\ast) < L(x,\gamma^\ast)$ .
    \label{saddle assumption}
\end{assumption}
 Assumption \ref{saddle assumption} is common in literature \cite{GOEBEL201716}.
\begin{assumption}
For every $x \in \mathcal{X}$, the function $L(x,\gamma) : \mathbb{R}^m \to \Bar{\mathbb{R}}$ is upper semicontinuous, and for every $\gamma \in \mathcal{Y}$, the function $L(x,\gamma) : \mathbb{R}^n \to \Bar{\mathbb{R}}$ is lower semicontinuous.
\label{lscusc}
\end{assumption}
Before proceeding to the next assumption, we define the (sub/super) differential $(\partial_x L(x,\gamma) \times \partial^+_\gamma L(x,\gamma)) $ for each $(x,\gamma) \in \mathcal{X} \times \mathcal{Y}$ as follows:
\begin{equation*}
    \partial_x L(x,\gamma) = \{ g_x \in \mathbb{R}^n \; | L(y,\gamma) \geq L(x,\gamma) + \inprod{g_x}{y-x} \; \; \forall \; y \in \mathbb{R}^n \}
\end{equation*}
and 
\begin{equation*}
    \partial_\gamma^+ L(x,\gamma) = \{ g_\gamma \in \mathbb{R}^m \; | L(x,\delta) \leq L(x,\gamma) + \inprod{g_\gamma}{\delta-\gamma} \; \forall \; \delta \in \mathbb{R}^m \}.
\end{equation*}
\begin{assumption}
  For each $(x,\gamma) \in \mathcal{X} \times \mathcal{Y}$, $\partial_xL(x,\gamma) \times \partial^+\gamma L(x,\gamma)$ is nonempty, and there exists a constant $G > 0$ such that $\forall \; g_x \in \partial_x L(x,\gamma)$, $\norm{g_x} \leq G$ and $\forall \; g_\gamma \in \partial_\gamma^+ L(x,\gamma)$,   $\norm{g_\gamma} \leq G$ .
    \label{subg}
\end{assumption}
The remainder of the paper consists of two main parts.

In first part, we develop a projected gradient dynamics on a Riemannian manifold to address problem \eqref{primal} and establish an equivalence between this projected dynamical system and the continuous-time saddle point mirror descent dynamics. Furthermore, we demonstrate that the solution of the projected dynamical system converges to $\mathcal{S}$ for any initial condition.

In the second part, we present stochastic first-order and zeroth-order saddle point mirror descent algorithms and utilize continuous-time analysis to prove the almost sure convergence of the iterates these algorithms.

In the next section, we develop projected dynamical  systems on Riemannian manifold to solve the saddle point problem \eqref{primal}.
\section{Projected Dynamical System in Non-Euclidean domain}
Before presenting the projected dynamical system aimed at solving problem  \eqref{primal}, it is essential to revisit the notion of a Riemannian metric for the sets $\mathcal{X}$ and $\mathcal{Y}$. 
We start by considering two $\sigma_R$-strongly convex functions, $R_x: \mathbb{R}^n \to \Bar{R}$ and $R_\gamma: \mathbb{R}^m \to \Bar{R}$. Given a set $\mathcal{X} \subseteq \mathbb{R}^n$, a Riemannian metric is a mapping that associates, to each point $x \in \mathcal{X}$, an inner product on the set $\mathbb{R}^n$ (more specifically, on $T_x \mathbb{R}^n$, the usual tangent space to $\mathbb{R}^n$ at $x$). The Riemannian metrics for each $x \in \mathcal{X}$ and $\gamma \in \mathcal{Y}$ on $\mathbb{R}^n$ and $\mathbb{R}^m$ are defined as follows:
\begin{equation}
\langle u_1, u_2 \rangle_x = u_1^\top \nabla^2 R_x(x) u_2 \quad \forall \; u_1,u_2 \in \mathbb{R}^n, \quad \text{and} \quad \langle v_1, v_2 \rangle_\gamma = v_1^\top \nabla^2 R_\gamma(\gamma) v_2 \quad \forall \; v_1,v_2 \in \mathbb{R}^m.
\label{inner}
\end{equation}
It can be verified that for each $x \in \mathcal{X}$ and $\gamma \in \mathcal{Y}$, this indeed defines a valid inner product on $\mathbb{R}^n$ and $\mathbb{R}^m$, respectively. Thus, for each $x \in \mathcal{X}$ and $\gamma \in \mathcal{Y}$, the inner product defined in \eqref{inner} induces a valid norm on $\mathbb{R}^n$ and $\mathbb{R}^m$, respectively. In other words, for all $x \in \mathcal{X}$ and $\gamma \in \mathcal{Y}$, we define the norms on the set $\mathbb{R}^n$ and $\mathbb{R}^m$ as:
\begin{equation*}
\norm{u}_x^2 = u^\top \nabla^2 R_x(x) u \quad \text{and} \quad \norm{v}_\gamma^2 = v^\top \nabla^2 R_\gamma(\gamma) v \; \; \; (u,v) \in \mathbb{R}^n \times \mathbb{R}^m.
\end{equation*}

Therefore, for each $x \in \mathcal{X}$, we obtain a finite-dimensional Hilbert space as $(T_x(\mathbb{R}^n)  = \mathbb{R}^n, \norm{.}_x )$. Similarly, we obtain one for each $\gamma \in \mathcal{Y}$.
Now, let's focus on the function $L(x,\gamma) : \mathbb{R}^n \times \mathbb{R}^m \to \Bar{\mathbb{R}}$ and define the Riemannian sub and super gradient of this function with respect to both $x$ and $\gamma$ utilizing the Riemannian metric specified in \eqref{inner}.
\\ We term $g_x' \in \mathbb{R}^n$ as a Riemannian subgradient of the function $L(x,\gamma)$ with respect to  $x \in \mathcal{X}$ if and only if
\begin{equation*}
\langle v, g_x' \rangle_x \leq D_x L(x,\gamma)[v] \quad \forall v \; \in \mathbb{R}^n.
\end{equation*}
Here, $D_x L(x,\gamma)[v]$ represents a directional derivative of the function $L$ with respect to $x$ in the direction $v$.
Similarly, we can define the Riemannian supergradient of the function $L(x,\gamma)$, denoted as $g_\gamma' \in \mathbb{R}^m$, with respect to $\gamma \in \mathcal{Y}$  if and only if:
\begin{equation*}
\langle w, g_\gamma' \rangle_\gamma \geq D_\gamma L(x,\gamma)[w] \quad \forall \; w \in \mathbb{R}^m.
\end{equation*}
It is straightforward to verify that $g_x' \in \mathbb{R}^n$ is a Riemannian subgradient of $L(x,\gamma)$ with respect to $x$ if and only if there exists $g_x \in \partial_x L(x,\gamma)$ such that $g_x' = \nabla^2 R_x(x)^{-1} g_x$. Similarly, $g_\gamma' \in \mathbb{R}^m$ is a supergradient of $L(x,\gamma)$ with respect to $\gamma$ if and only if there exists $g_\gamma \in \partial_\gamma L(x,\gamma)$ such that $g_\gamma' = \nabla^2 R_\gamma(\gamma)^{-1} g_\gamma$.
\\ Define the sets
\begin{equation*}
\partial_x' L(x,\gamma) = \{ g_x' \mid g_x' \text{ is a Riemannian subgradient} \}
\end{equation*}
and
\begin{equation*}
\partial_\gamma'^{+} L(x,\gamma) = \{ g_\gamma' \mid g_\gamma' \text{ is a Riemannian supergradient} \}.
\end{equation*}
From the preceding discussion, we conclude that  $\partial_x' L(x,\gamma) = \nabla^2R_x(x)^{-1} \partial_x L(x,\gamma)$, and likewise, \\ $\partial_\gamma'^{+} L(x,\gamma) = \nabla^2R_\gamma(\gamma)^{-1} \partial_\gamma^+ L(x,\gamma)$.

In this context, we define the normal cone of the set $\mathcal{X}$ at point $x$ $\in$ $\mathcal{X}$ (denoted as $\mathcal{N}_\mathcal{X}'(x)$), considering the Riemannian metric, as follows:
\begin{equation*}
    \mathcal{N}_\mathcal{X}'(x) = \{\eta_x' \in \mathbb{R}^n \; | \inprod{\eta_x'}{\nu_x}_x \leq 0 \; , \; \; \forall \; \nu_x \in \mathcal{T}_\mathcal{X}(x)\} 
\end{equation*}

Clearly, in this definition, the normal cone represents the polar cone of the tangent cone with respect to the Riemannian metric.  Similarly, we  define the normal cone of the set $\mathcal{Y}$ at point $\gamma$ $\in$ $\mathcal{Y}$ (denoted as $\mathcal{N}_\mathcal{Y}'(\gamma)$) under the Riemannian metric as follows:
\begin{equation*}
    \mathcal{N}_\mathcal{Y}'(\gamma) = \{\eta_\gamma' \in \mathbb{R}^m \; | \inprod{\eta_\gamma'}{\nu_\gamma}_\gamma \leq 0 \; , \; \; \forall \; \nu_\gamma \in \mathcal{T}_\mathcal{Y}(\gamma) \} .
\end{equation*}

We can straightforwardly verify that 
$ \mathcal{N}_\mathcal{X}'(x) = \nabla^2 R_x(x)^{-1} \mathcal{N}_\mathcal{X}(x)$ and $\mathcal{N}_\mathcal{Y}'(\gamma) = \nabla^2 R_\gamma(\gamma)^{-1} \mathcal{N}_\mathcal{Y}(\gamma)$. 

Based on the discussion thus far in this section, we are now prepared to present the projected saddle point dynamics within the non-Euclidean domain to address the problem given  in \eqref{primal}. The formulation of these dynamics are as follows:

\begin{equation}
\begin{split}
 &   \Dot{x}(t) \in \mathcal{P}_{\mathcal{T}_\mathcal{X}(x(t))}^{x(t)}(-\partial_x' L(x(t),\gamma(t)) ) = \{ \nu_1 \in \mathbb{R}^n \; | \;  \nu_1 = \argmin\limits_{\nu \;   \in \mathcal{T}_\mathcal{X}(x(t))} \norm{\nu + g_x'(t)}_{x(t)}^2  ,\;  g_x'(t) \in \partial_x' L(x(t),\gamma(t)) \}
    \end{split}
    \label{projected primal}
\end{equation}
\begin{equation*}
      x(0) = x_0 \in \mathcal{X}.
\end{equation*}
\begin{equation}
    \Dot{\gamma}(t) \in \mathcal{P}_{\mathcal{T}_\mathcal{Y}(\gamma(t))}^{\gamma(t)}(\partial_\gamma'^{+} L(x(t),\gamma(t))) = \{ \nu_2 \in \mathbb{R}^m \; | \;  \nu_2 = \argmin\limits_{\nu \in \mathcal{T}_\mathcal{Y}(\gamma(t))} \norm{\nu - g_\gamma'(t)}_{\gamma(t)}^2  , \; g_\gamma'(t) \in \partial_\gamma'^{+} L(x(t),\gamma(t)) \}
    \label{projected dual}
\end{equation}
\begin{equation*}
    \gamma(0) = \gamma_0 \in \mathcal{Y}.
\end{equation*}
Note that the projections onto the sets $\mathcal{T}_\mathcal{X}(x(t))$ and $\mathcal{T}_\mathcal{Y} (\gamma(t))$ are explicitly formulated with consideration for the metric induced by the local norm associated with $(x(t),\gamma(t)) \in \mathcal{X} \times \mathcal{Y}$, as defined in \eqref{inner}. The following Proposition extends Moreau's Theorem (as presented in Theorem $3.2.5$ of \cite{hiriart2013convex}) in the context of Riemannian Metrics. The proof can be conducted in a similar manner as outlined in \cite{hiriart2013convex}, with the adjustment of replacing the Euclidean metric with the Riemannian metric.
\begin{proposition}[\cite{hiriart2013convex}]
  Let $d \in \mathbb{R}^n$.   The following two statements are equivalent.
    \begin{enumerate}
        \item $\nu = \mathcal{P}_{\mathcal{T}_\mathcal{X}(x)}^x(d)$ and $\eta = \mathcal{P}_{\mathcal{N}_\mathcal{X}'(x)}^x(d) $.
        \item $\nu + \eta = d$, where $\nu \in \mathcal{T}_\mathcal{X}(x)$ and $ \eta \in \mathcal{N}_\mathcal{X}'(x) $, and $\inprod{\eta}{\nu}_x = 0$.
    \end{enumerate}
    \label{Morea}
\end{proposition}

\begin{remark}
 The dynamics described in \eqref{projected primal}-\eqref{projected dual} constitute a coupled  projected differential inclusions. However, a critical concern arises regarding the well-defined nature of these ordinary differential inclusions. This stems from observations indicating that the tangent cones $\mathcal{T}_{\mathcal{X}}(x)$ and $\mathcal{T}_{\mathcal{Y}}(\gamma)$ are not well-defined if $(x, \gamma) \notin \mathcal{X}\times \mathcal{Y}$. Nevertheless, these discrepancies can be resolved if we can establish the existence of a viable maximal Caratheodory solution for the dynamics associated with \eqref{projected primal}-\eqref{projected dual}. In other words, we have to show that $\forall \; (x_0,\gamma_0) \in \mathcal{X} \times \mathcal{Y}$, there exists absolutely continuous function $(x(t),\gamma(t)) : [0,\infty) \to \mathcal{X} \times \mathcal{Y}$ which satisfies \eqref{projected primal}-\eqref{projected dual} a.e. In the subsequent section we establish that.
 \end{remark}
 \begin{remark}
It is noteworthy that the dynamics in \eqref{projected primal} and \eqref{projected dual} exhibit a nuanced level of generality when compared to the formulations in \cite{doi:10.1137/18M1229225}. This distinction arises primarily due to the presence of set-valued maps in these formulations. This deviation renders the existence theorem provided in \cite{doi:10.1137/18M1229225} not directly applicable in this context. However, drawing inspiration from the analytical framework presented in that paper, we present the following lemma, providing an equivalent representation of the dynamics in terms of the normal cone rather than the tangent cone. This lemma plays a crucial role in establishing the existence of Caratheodory solutions for the projected dynamical system outlined in \eqref{projected primal} and \eqref{projected dual}.   The proof of this lemma follows procedures akin to those detailed in \cite{doi:10.1137/18M1229225}. For clarity and completeness, we provide a step-by-step account of the proof.
\end{remark}
\subsection{An Equivalent Representation}
\begin{lemma}
    Consider the ordinary differential inclusions given as follows:
    \begin{equation}
\begin{split}
  &  \Dot{x}(t) \in  - (  \partial_x' L(x (t), \gamma(t)) + \mathcal{N}_\mathcal{X}'(x(t))) =  - \nabla ^2  R_x 
(x(t))^{-1}    (  \partial_x L(x (t), \gamma(t)) + \mathcal{N}_\mathcal{X}(x(t)))
\\ & x(0) = x_0 \in \mathcal{X}
    \label{primal dynamics}
    \end{split}
\end{equation}
and 
\begin{equation}
    \begin{split}
&        \Dot{\gamma} (t) \in (  \partial'^+_\gamma L(x(t), \gamma(t)) -\mathcal{N}_\mathcal{Y}'(\gamma(t)))=  
\nabla ^2  R_\gamma
(\gamma(t))^{-1}    (  \partial^+_\gamma L(x(t), \gamma(t)) -\mathcal{N}_\mathcal{Y}(\gamma(t)))
\\ & \gamma(0) = \gamma_0 \in \mathcal{Y}.
    \end{split}
    \label{dual dynamics}
\end{equation} 
   The viable  Caratheodory solutions of \eqref{projected primal}-\eqref{projected dual}, if they exist, for any initial condition $(x_0, \gamma_0) \in \mathcal{X} \times \mathcal{Y}$, will also be Caratheodory solutions of the coupled ordinary differential inclusions presented in \eqref{primal dynamics}-\eqref{dual dynamics}, and vice versa.
    \label{equivalence}
\end{lemma}
 \begin{proof}
        Suppose $(x(t),\gamma(t))$ is a viable Caratheodory solution of \eqref{projected primal}-\eqref{projected dual}. In other words,
    $\exists \; g_x'(t) \in \partial_x' L(x(t),\gamma(t))$  and $g_\gamma'(t) \in \partial_\gamma'^+  L(x(t),\gamma(t))$ for almost all $t$ such that
    \begin{equation}
        \Dot{x}(t) = \mathcal{P}_{\mathcal{T}_\mathcal{X}(x(t))}^x(-g_x'(t)) \; \; \text{and} \; \; \Dot{\gamma}(t)  = \mathcal{P}_{\mathcal{T}_\mathcal{Y}(\gamma(t))}^\gamma(g_\gamma'(t)).
        \label{1}
    \end{equation}
    
 According to Proposition \ref{Morea},  $\exists \; \eta_x'(t) \in \mathcal{N}_{\mathcal{X}}' (x(t))$ and $\exists \; \eta_\gamma'(t) \in \mathcal{N}_{\mathcal{Y}}' (\gamma(t))$ for almost all $t$, such that:
 \begin{equation*}
     \Dot{x}(t) = - g_x'(t) - \eta_x'(t) \; \; \text{and} \; \;  \Dot{\gamma}(t) =  g_\gamma'(t) - \eta_\gamma'(t).
 \end{equation*}

    Thus, $(x(t),\gamma(t))$ is also a Caratheodory solution of ordinary differential inclusions \eqref{primal dynamics}-\eqref{dual dynamics}.
    Conversely, assume that $(x(t) , \gamma(t))$ is a viable  Caratheodory solution of ordinary differential inclusions \eqref{primal dynamics}-\eqref{dual dynamics}.
    In that case, for almost all $t \in [0,\infty)$ $\exists \;  g_x(t) \in \partial_x L(x(t),\gamma(t))$ and $\exists \; g_\gamma(t) \in \partial^+_\gamma L(x(t),\gamma(t))$ and $\eta_x(t) \in \mathcal{N}_\mathcal{X}(x(t))$ , $\eta_\gamma(t) \in \mathcal{N}_\mathcal{Y}(\gamma(t))$ such that 
    \begin{equation*}
        \Dot{x}(t) =   -\nabla^2 R_x(x(t))^{-1} (g_x(t)  + \eta_x(t)). 
    \end{equation*}
    and 
    \begin{equation*}
         \Dot{\gamma}(t) =  \nabla^2 R_\gamma(\gamma(t))^{-1}( g_\gamma(t) -\eta_\gamma(t)).
    \end{equation*}
     We next show that \eqref{1} hold. However, \eqref{1} holds if and only if 
     \begin{equation}
        \inprod{\nabla^2 R_x(x(t)) (\Dot{x}(t)+\nabla^2 R_x(x(t))^{-1} g_x(t))}{v -\Dot{x}(t)} \geq 0 \; \; \forall \; \; v \in \mathcal{T}_\mathcal{X}(x(t)).
        \label{rumi1}
    \end{equation}
    Consider the LHS of \eqref{rumi1} and plug $\Dot{x}(t) =   -\nabla^2 R_x(x(t))^{-1} (g_x(t)  + \eta_x(t))$ resulting in
    \begin{equation}
    \begin{split}
         &  \inprod{\nabla^2 R_x(x(t)) (\Dot{x}(t)+\nabla^2 R_x(x(t))^{-1} g_x(t))}{v -\Dot{x}(t)}
          =  \inprod{-\eta_x(t)}{v} + \inprod{\Dot{x}(t)}{\eta_x(t)} \geq 0 .
          \end{split}
          \label{mimia}
    \end{equation}
    The inequality in \eqref{mimia} holds true owing to the definition of the normal cone and the fact that $-\dot{x}(t) \in \mathcal{T}_\mathcal{X}(x(t))$, which can be verified through the following argument:
\\ Since $x(t)$ is absolutely continuous hence it is differentiable almost everywhere . Hence for almost all $t$ we have
\begin{equation*}
    \Dot{x}(t) = \lim\limits_{\Delta t \downarrow 0} \frac{x(t)-x(t-\Delta t)}{\Delta t}.
\end{equation*}
Choose any sequence $t_k \downarrow 0$ we have 
\begin{equation*}
    \Dot{x}(t) = - \lim\limits_{t_k \downarrow 0} \frac{x(t-t_k)-x(t)}{t_k} \in - \mathcal{T}_\mathcal{X}(x(t)).
\end{equation*}
Hence \eqref{mimia} holds. Similarly, we can prove for $\gamma$. So, $(x(t),\gamma(t))$ is a Caratheodory solution of the projected dynamical system \eqref{projected primal}-\eqref{projected dual}.
\end{proof}
The equivalency established by lemma \ref{equivalence} provides an equivalent viewpoint on the projected differential inclusions in \eqref{projected primal}-\eqref{projected dual}.
This equivalence draws parallels with the continuous-time mirror descent dynamics introduced in \cite{10.1145/3188745.3188798}, designed for online bandit problems featuring continuously differentiable objective functions. However, our study diverges significantly as it delves into the stability analysis of saddle point dynamics,  and notably, we refrain from imposing any assumptions regarding the differentiability of the objective function. In light of these considerations, we coin the dynamics presented in \eqref{primal dynamics}-\eqref{dual dynamics} as continuous-time saddle point mirror descent dynamics. This terminology is apt, acknowledging the influence of prior works while highlighting the distinctive features of our approach. These discussions unveil a projected dynamical viewpoint for continuous-time saddle point mirror descent dynamics in \eqref{primal dynamics}-\eqref{dual dynamics}, offering a more intuitive understanding in geometric terms.
\\ However, concerns regarding the well-posed nature of the associated ordinary differential inclusions \eqref{primal dynamics}-\eqref{dual dynamics} persist. In addressing this concern, the subsequent lemma not only demonstrates the existence of viable Caratheodory solutions for the ordinary differential inclusions \eqref{primal dynamics}-\eqref{dual dynamics} but also alleviates any apprehensions about the well-posedness of the dynamics outlined in \eqref{primal dynamics}-\eqref{dual dynamics}. This, in turn, establishes the well-posed nature of the projected ordinary differential inclusions in \eqref{projected primal}-\eqref{projected dual}.

\subsection{Existence of Viable Caratheodory Solutions of \eqref{primal dynamics}-\eqref{dual dynamics}}

Before delving into the main existence theorem, we make the following observation.

For every $x \in \mathcal{X}$, we have defined a Hilbert space as $(\mathbb{R}^n,\norm{.}_x)$. Since norms on finite-dimensional spaces are equivalent, there exist four functions: $\kappa_1,\kappa_2 : \mathcal{X} \to \mathbb{R}^{+}$ and $\kappa_3,\kappa_4: \mathcal{Y} \to \mathbb{R}^{+}$ leading to  the following relation: 
\begin{equation*}
    \kappa_2 (x) \norm{u_1}_2 \leq \norm{u_1}_x \leq \kappa_1 (x) \norm{u_1}_2 \; \forall \; u _1 \in \mathbb{R}^n \; \text{and}  \; \; \kappa_4 (\gamma) \norm{u_2}_2 \leq \norm{u_2}_\gamma \leq \kappa_3(\gamma)  \norm{u_2}_2 \; \forall \; u_2 \in \mathbb{R}^m.
\end{equation*}
In the upcoming lemma, we'll demonstrate that these functions $\kappa_1$, $\kappa_2$, $\kappa_3$ and $\kappa_4$ are independent of both $x$ and $\gamma$, given the compactness of $\mathcal{X}$ and $\mathcal{Y}$. 
\begin{lemma}
    There exists $\kappa_1 , \kappa_2 > 0$ such that for all $x \in \mathcal{X}$ and $\gamma \in \mathcal{Y}$ we have 
    \begin{equation*}
        \kappa_2 \norm{u_1}_2 \leq \norm{u_1}_x \leq \kappa_1 \norm{u_1}_2 \; \forall \; u_1 \in \mathbb{R}^n \; \text{and} \; \;  \kappa_2 \norm{u_2}_2 \leq \norm{u_2}_
        \gamma \leq \kappa_1 \norm{u_2}_2 \; \forall \; u_2 \in \mathbb{R}^m.
    \end{equation*}
    \label{normequivalence}
\end{lemma}
\begin{proof}
    Recall that the maximum eigenvalue of $\nabla^2 R_x(x)$ is given by:

\begin{equation*}
    \lambda_{\max} (\nabla^2 R_x(x)) = \max_{\lVert u \rVert_2 = 1} u^\top \nabla^2 R_x(x) u.
\end{equation*}

Now, consider any $u \in \mathbb{R}^n$, then:

\begin{equation*}
    \lambda_{\max} (\nabla^2 R_x(x)) \geq \left(\frac{u}{\lVert u \rVert}\right)^\top \nabla^2 R_x(x) \left(\frac{u}{\lVert u \rVert}\right).
\end{equation*}

In other words:

\begin{equation*}
    \lVert u \rVert_x^2 \leq \lambda_{\max} (\nabla^2 R_x(x)) \lVert u \rVert_2^2.
\end{equation*}

Choose $\kappa_1 = \max\limits_{x \in \mathcal{X}} (\lambda_{\max} (\nabla^2 R_x(x)))$. Therefore, there exists $\kappa_1 > 0$ such that:

\begin{equation*}
    \lVert u \rVert_x^2 \leq \kappa_1 \lVert u \rVert_2^2 \quad \forall \; u \in \mathbb{R}^n.
\end{equation*}

Without loss of generality, we can assume that:

\begin{equation*}
\lVert v \rVert_\gamma^2 \leq \kappa_1 \lVert v \rVert_2^2 \quad \forall \; v \in \mathbb{R}^m.
\end{equation*}

Similarly, we can prove that there exists $\kappa_2 > 0$ such that:

\begin{equation*}
\lVert u \rVert_x^2 \geq \kappa_2 \lVert u \rVert_2^2 \quad \forall \; u \in \mathbb{R}^n \quad \text{and} \quad \lVert v \rVert_\gamma^2 \geq \kappa_2 \lVert v \rVert_2^2 \quad \forall \; v \in \mathbb{R}^m.
\end{equation*}
\end{proof}

In the next lemma, we prove the existence of Caratheodory solution of the ordinary differential inclusions given in \eqref{primal dynamics} and \eqref{dual dynamics}.

\begin{lemma} 
For all initial condition $(x_0,\gamma_0) \in \mathcal{X} \times \mathcal{Y}$,
there exists an absolutely continuous map $t \mapsto (x(t),\gamma(t)) \in  \mathcal{X}\times \mathcal{Y}$ satisfying \eqref{primal dynamics}- \eqref{dual dynamics} a.e. 
    \label{existence}
 \end{lemma}
 \begin{proof}
     Consider the differential inclusions in coupled form as
     \begin{equation}
         \begin{bmatrix}
             \Dot{x}
             \\ \Dot{\gamma} 
         \end{bmatrix} \in \begin{bmatrix}
         -     \nabla ^2  R_x 
(x)^{-1}    (  \partial_x L(x, \gamma) + \mathcal{N}_\mathcal{X}(x)) \\ 
\nabla ^2  R_\gamma
(\gamma)^{-1}    (  \partial^{+}_\gamma L(x, \gamma) -\mathcal{N}_\mathcal{Y}(\gamma))
         \end{bmatrix}.
         \label{ext}
     \end{equation}

 Now we are ready to show the existence of Caratheodory solutions of \eqref{ext}. 
 Instead of directly establishing the existence of a viable Caratheodory solution for \eqref{ext}, we introduce the following ordinary differential inclusion:
     \begin{equation}
         \begin{bmatrix}
             \Dot{x}
             \\ \Dot{\gamma} 
         \end{bmatrix} \in \begin{bmatrix}
         -     \nabla ^2  R_x 
(x)^{-1}    (  \partial_x L(x, \gamma) + \Hat{\mathcal{N}}_\mathcal{X}(x)) \\ 
\nabla ^2  R_\gamma
(\gamma)^{-1}    (  \partial^{+}_\gamma L(x, \gamma) -\Hat{\mathcal{N}}_\mathcal{Y}(\gamma))
         \end{bmatrix}. 
         \label{ext1}
         \end{equation}
         where,
         \begin{equation*}
             \begin{split}
                 &\Hat{\mathcal{N}}_\mathcal{X}(x) = \{ \eta_x \in \mathcal{N}_\mathcal{X}(x) \; | \; \norm{\eta_x} \leq \frac{\kappa_1}{\kappa_2}G \} 
                 \\ & \Hat{\mathcal{N}}_\mathcal{Y}(\gamma) = \{ \eta_\gamma \in \mathcal{N}_\mathcal{Y}(\gamma) \; | \; \norm{\eta_\gamma} \leq  \frac{\kappa_1}{\kappa_2}G\}.
             \end{split}
         \end{equation*}
      It is evident that a Caratheodory solution of \eqref{ext1} also be a Caratheodory solution for \eqref{ext}. To enhance clarity, denote right-hand side of equation \eqref{ext1} as $F(x,\gamma)$. Our objective is to establish the hypothesis  outlined in Proposition \ref{dh} for $F(x,\gamma)$, which in turn implies the existence of a Caratheodory solution for \eqref{ext1}.
\\ \textbf{Claim-1:} $F(x,\gamma)$ is an upper semi-continuous (usc) map for each $(x,\gamma) \in \mathcal{X} \times \mathcal{Y}$.

Assume that $F$ is not usc at $(x,\gamma) \in \mathcal{X} \times \mathcal{Y}$. Then there exists $\epsilon > 0$ such that $\forall \; \delta > 0$: $\exists \; (x',\gamma') \in U((x,\gamma),\delta) $ and $\exists \; (\mu',\nu') \in F(x',\gamma')$ where $(\mu',\nu') \notin F(x,\gamma) + U(0,2\epsilon)$. In other words 
\begin{equation*}
    \norm{(\mu',\nu')-(\mu,\nu)} \geq \epsilon \; \; \forall \; (\mu,\nu) \in F(x,\gamma).
\end{equation*}
Choose $\delta = \frac{1}{n}$, $(n \in \mathbb{N})$ then $\exists \; (x_n,\gamma_n) \in U((x,\gamma),
\frac{1}{n}) $ and $\exists \; (\mu_n,\nu_n) \in F(x_n,\gamma_n)$ 
\begin{equation}
    \norm{(\mu_n,\nu_n)-(\mu,\nu)} \geq \epsilon \; \; \forall \; (\mu,\nu) \in F(x,\gamma) \; \; \forall \; n \in \mathbb{N}.
    \label{contra}
\end{equation}
Clearly, $(x_n,\gamma_n) \to (x,\gamma)$. From the definition of $F$ we have 
\begin{equation*}
    \begin{split}
&   \mu_n \in   -   \nabla ^2  R_x 
(x_n)^{-1}    (  \partial_x L(x_n, \gamma_n) + \Hat{\mathcal{N}}_\mathcal{X}(x_n))
\\ & \nu_n \in \nabla ^2  R_\gamma
(\gamma_n)^{-1}    (  \partial^{+}_\gamma L(x_n, \gamma_n) -\Hat{\mathcal{N}}_\mathcal{Y}(\gamma_n)).
    \end{split}
\end{equation*}
That implies $\exists \; g_x(n) \in \partial_x L(x_n,\gamma_n)$ and $\eta_x(n) \in \Hat{\mathcal{N}}_\mathcal{X}(x_n)$ such that
\begin{equation*}
    \mu_n = - \nabla^2 R_x(x_n)^{-1} (g_x(n)+ \eta_x(n)).
\end{equation*}
Similarly,
\begin{equation*}
    \nu_n = \nabla^2 R_\gamma(\gamma_n)^{-1} (g_\gamma(n) - \eta_\gamma(n)).
\end{equation*}
where, $g_\gamma(n) \in \partial^{+}_\gamma L(x_n,\gamma_n)$ and $\eta_\gamma(n) \in \Hat{\mathcal{N}}_\mathcal{Y}(\gamma_n)$. Since $g_x(n)$ and $\eta_x(n)$  are bounded, we can safely assume that $g_x(n) + \eta_x(n) \to \Bar{g}_x+\Bar{\eta}_x$ (as otherwise there exists a convergent subsequence). Hence, $\{\mu_n\} \to \Bar{\mu}$ where,
\begin{equation}
    \Bar{\mu} = -\nabla^2R_x(x)^{-1} (\Bar{g}_x+\Bar{\eta}_x).
    \label{151}
\end{equation}
Similarly, $\nu_n \to \Bar{\nu}$ such that 
\begin{equation}
    \Bar{\nu} = \nabla^2 R_\gamma(\gamma)^{-1} (\Bar{g}_\gamma-\Bar{\eta}_\gamma)
    \label{161}
\end{equation}
where, $g_\gamma(n) - \eta_\gamma(n) \to \Bar{g}_\gamma-\Bar{\eta}_\gamma$ for a similar reason already mentioned.
Since $g_x(n) \in \partial_x L(x_n,\gamma_n)$ from the definition of subgradient
\begin{equation}
    L(y,\gamma_n) \geq L(x_n,\gamma_n) + \inprod{g_x(n)}{y-x_n} \; \; \forall \; y \in \mathcal{X}
    \label{15}
\end{equation}
and from the definition of supergradient we get
\begin{equation}
    -L(x_n,\gamma) \geq - L(x_n,\gamma_n) +\inprod{-g_\gamma(n)}{\gamma-\gamma_n}.
    \label{suos}
\end{equation}
Hence, from \eqref{15} and \eqref{suos} we have
\begin{equation}
    L(y,\gamma_n) \geq L(x_n,\gamma) + \inprod{g_x(n)}{y-x_n} - \inprod{g_\gamma(n)}{\gamma-\gamma_n}.
    \label{16}
\end{equation}
Letting $n \to \infty $ on both sides of \eqref{16} we get
\begin{equation}
    \limsup\limits_{n\to \infty} L(y,\gamma_n) \geq \liminf\limits_{n \to \infty} L(x_n,\gamma) + \inprod{\Bar{g}_x}{y-x}.
    \label{17}
\end{equation}
Since $L$ is upper semicontinuous with respect to the variable $\gamma$ and hence from the definition of upper semicontinuity $L(y, \gamma) \geq \limsup\limits_{n \to \infty} L(y,\gamma_n)$. Similarly from the definition of lower semicontinuity we get $\liminf\limits_{n \to \infty} L(x_n,\gamma)  \geq L(x,\gamma) $. So, from \eqref{17}, we obtain 
\begin{equation*}
    L(y,\gamma) \geq L(x,\gamma)+ \inprod{\Bar{g}_x}{y-x}.
\end{equation*}
Hence $\Bar{g}_x \in \partial_x L(x,\gamma)$. In a similar fashion we can show that $\Bar{g}_\gamma \in \partial^+_\gamma L(x,\gamma)$.
\\ Now from the definition of normal cone 
\begin{equation*}
    \begin{split}
        & \inprod{\eta_x(n)}{y-x_n} \leq 0 \; \; \; \forall \; y \in \mathcal{X}
        \\ &  \inprod{\Bar{\eta}_x}{y-x} \leq 0 \; \; \;  \; \; \; \; \; \text{since} \; \; \eta_x(n) \to \Bar{\eta}_x.
    \end{split}
\end{equation*}
Also $\norm{\Bar{\eta}_x} \leq  \frac{\kappa_1}{\kappa_2}G$. Hence $\Bar{\eta}_x \in \Hat{\mathcal{N}}_\mathcal{X}(x)$. Similarly, $\Bar{\eta}_\gamma \in \Hat{\mathcal{N}}_\mathcal{Y}(\gamma)$. 
Hence from \eqref{151} and \eqref{161} we get that $(\Bar{\mu},\Bar{\nu}) \in F(x,\gamma)$. However, this is a contradiction to \eqref{contra} considering the fact $(\mu_n,\nu_n) \to (\Bar{\mu},\Bar{\nu})$. Hence, $F(x,\gamma)$ is upper semicontinuous $\forall \; (x,\gamma) \in \mathcal{X}\times \mathcal{Y}$.
\\ \textbf{Claim-2:} For each $(x,\gamma) \in \mathcal{X} \times \mathcal{Y}$, $F(x,\gamma)$ is compact and convex. 
\\ For convexity, it is enough to show that both $\Hat{\mathcal{N}}_\mathcal{X}(x)$ and  $\Hat{\mathcal{N}}_\mathcal{Y}(\gamma)$ are convex  for each $x \in \mathcal{X}$ , $\gamma \in \mathcal{Y}$ considering the fact sum of two convex sets is convex. 
\\ Suppose $\eta_1 , \eta_2 \in $ $\Hat{\mathcal{N}}_\mathcal{X}(x) $ then for any $\theta \in (0,1)$ we have
\begin{equation*}
    \theta \eta_1 + (1-\theta) \eta_2 \in \mathcal{N}_\mathcal{X} (x) \; \; \; \; \text{and} \; \;    \norm{\theta \eta_1 + (1-\theta) \eta_2} \leq \frac{\kappa_1}{\kappa_2}G
\end{equation*}
which implies $\theta \eta_1 + (1-\theta) \eta_2 \in \Hat{\mathcal{N}}_\mathcal{X}(x)$. In a similar way, we can show that $\Hat{\mathcal{N}}_\mathcal{Y}(\gamma)$ is also convex. 
\\ Since $\forall \; x \in \mathcal{X}$ and $\forall \; \gamma \in
 \mathcal{Y}$ both the sets $\partial_x L(x,y)$ and $\Hat{\mathcal{N}}_\mathcal{X}(x)$ are compact, their sum is  also  compact. Hence, for each $(x,\gamma) \in \mathcal{X} \times \mathcal{Y}$, $F(x,\gamma)$ is convex and compact. 
 \\ \textbf{Claim-3:} For each $(x,\gamma) \in \mathcal{X} \times \mathcal{Y}$,  we have
 \begin{equation*}
     \sup\limits_{y \in F(x,\gamma)} \norm{y} \leq (1 + \frac{\kappa_1}{\kappa_2})G.
 \end{equation*}
 This claim is true in view of Assumption \ref{subg} and from the definition of the sets $\Hat{\mathcal{N}}_\mathcal{X}(x)$ and $\Hat{\mathcal{N}}_\mathcal{Y}(\gamma)$.
\\ \textbf{Claim-4:} $\forall \; (x,\gamma) \in \mathcal{X} \times \mathcal{Y}$
\begin{equation*}
    F(x,\gamma) \cap \mathcal{T}_{\mathcal{X}\times \mathcal{Y}}(x,\gamma) \neq \emptyset.
\end{equation*}
It is enough to prove the claim  for $\forall \; (x,\gamma) \in bd(\mathcal{X}\times \mathcal{Y})$. Consider $g_x' \in \partial_x' L(x,\gamma)$ and $g_\gamma' \in \partial_\gamma'^+ L(x,\gamma)$ and the following:
\begin{equation}
    v' = \mathcal{P}_{\mathcal{T}_\mathcal{X}(x)}^x (-g_x') 
    \label{opt1}
\end{equation}
and 
\begin{equation}
    w' = \mathcal{P}_{\mathcal{T}_\mathcal{Y}(\gamma)}^\gamma (g_\gamma').
    \label{opt2}
\end{equation}
Applying the optimality condition to \eqref{opt1}, we have that for all  $ y \in \mathcal{T}_\mathcal{X}(x)$
\begin{equation}
    \inprod{(v'+g_x')}{y-v'}_x \geq 0.
    \label{cone}
\end{equation}
Setting $y = 2 v' \in \mathcal{T}_\mathcal{X}(x)$ in \eqref{cone}, we obtain
\begin{equation}
       \inprod{(v'+g_x')}{v'}_x \geq 0.
      \label{cone1}
\end{equation}
Similarly, by setting $y = \frac{1}{2} v' \in \mathcal{T}_\mathcal{X}(x)$ in \eqref{cone},  we get  
\begin{equation}
     \inprod{(v'+g_x')}{v'}_x \leq 0.
      \label{cone2}
\end{equation}
Combining \eqref{cone1} and \eqref{cone2}, we deduce
\begin{equation*}
   \inprod{(v'+g_x')}{v'}_x = 0.
\end{equation*}
Thus, from \eqref{cone}, it follows that for all $ y \in \mathcal{T}_\mathcal{X}(x)$ 
\begin{equation*}
     \inprod{(v'+g_x')}{y}_x \geq 0.
\end{equation*}
In other words from the definition of normal cone in the Riemannian metric we have $(v'+g_x') \in - \mathcal{N}_\mathcal{X}'(x)$, that is, $\exists \; \eta_x' \in \mathcal{N}_\mathcal{X}'(x)$
\begin{equation*}
\begin{split}
    & v' + g_x' = - \eta_x'
     =>  v' = -\nabla^2 R_x(x)^{-1} (g_x+\eta_x) 
    \end{split}
\end{equation*}
where, $g_x \in \partial_x L(x,\gamma)$ and $\eta_x \in \mathcal{N}_\mathcal{X}(x)$.
Since $\mathcal{T}_\mathcal{X}(x)$ is a closed  cone, $0 \in \mathcal{T}_\mathcal{X}(x)$ and from \eqref{opt1} we get 
\begin{equation}
    \begin{split}
        & \norm{v'+g_x'}_x \leq \norm{g_x'}_x.
    \end{split}
    \label{ine}
\end{equation}
Putting $v' = -(g_x'+ \eta_x')$, we obtain: 
\begin{equation*}
    \begin{split}
        \norm{\eta_x'}_x  \leq \norm{g_x'}_x 
        = >     \norm{\eta_x} \leq  \frac{\kappa_1}{\kappa_2} G.
    \end{split}
\end{equation*}
Hence we show that $v' \in  -     \nabla ^2  R_x 
(x)^{-1}    (  \partial_x L(x, \gamma) + \Hat{\mathcal{N}}_\mathcal{X}(x))  \cap \mathcal{T}_\mathcal{X}(x) $. In a similar fashion we can demonstrate that $w' \in  \nabla ^2  R_\gamma
(\gamma)^{-1}    (  \partial^+_\gamma L(x, \gamma) -\Hat{\mathcal{N}}_\mathcal{Y}(\gamma)) \cap\mathcal{T}_\mathcal{Y}(\gamma)$.  Thus,  $F(x,\gamma) \cap \mathcal{T}_{\mathcal{X}\times \mathcal{Y}}(x,\gamma) \neq \emptyset$. Hence, from Proposition \ref{dh} we can conclude that for each $(x_0,\gamma_0) \in \mathcal{X} \times \mathcal{Y}$ there exists an absolutely continuous map $(x(t),\gamma(t)) : [0,\infty) \to \mathcal{X} \times \mathcal{Y}$ satisfying \eqref{primal dynamics}-\eqref{dual dynamics} a.e.
 \end{proof}
In the next subsection, we demonstrate that the set of saddle points serves as the equilibrium points for the continuous-time saddle point mirror descent dynamics \eqref{primal dynamics}-\eqref{dual dynamics}. Furthermore, we establish that all equilibrium points are Lyapunov stable, and all Caratheodory solutions of the differential inclusions \eqref{primal dynamics}-\eqref{dual dynamics} will converge to the set of saddle points $\mathcal{S}$.
\subsection{Stability of the Continuous Time Saddle Point Mirror Descent Dynamics} 
   \begin{lemma}
       The set of saddle points $\mathcal{S}$ coincides with the set of equilibrium points of the continuous time saddle point mirror descent dynamics given in \eqref{primal dynamics}-\eqref{dual dynamics}.
       \label{eq}
   \end{lemma}
   \begin{proof}
       The set of saddle points $\mathcal{S}$ can be characterized in the following way
       \begin{equation*}
         (x^\ast,\gamma^\ast) \in \mathcal{S} <=>  L(x^\ast,\gamma) \leq L(x^\ast, \gamma^\ast) \leq L(x,\gamma^\ast) \; \; \forall \; (x,\gamma) \in \mathcal{X} \times \mathcal{Y}.
       \end{equation*}
       Hence, 
       \begin{equation}
           x^\ast = \argmin\limits_{x \in \mathcal{X}} L(x,\gamma^\ast).
           \label{saddle1}
       \end{equation}
       Since $L(x,\gamma)$ is convex with respect to $x$, \eqref{saddle1} is true if and only if (by applying optimality condition) $\exists \; g_x^\ast \in \partial_x L(x^\ast,\gamma^\ast)$ such that 
       \begin{equation*}
       \begin{split}
          & \inprod{g_x^\ast}{x-x^\ast} \geq 0
          \\ <=> \; & 0 \in  \nabla^2 R_x(x)^{-1} (\partial_x L(x^\ast,\gamma^\ast) + \mathcal{N}_\mathcal{X}(x)).
           \end{split}
       \end{equation*}
       The last line is true because of strong convexity of $R_x$. On the other hand 
       \begin{equation}
          \gamma^\ast = \argmax\limits_{\gamma \in \mathcal{Y}} L(x^\ast,\gamma) .
          \label{saddle2}
       \end{equation}
       Since $-L$ is convex with respect to $\gamma$, \eqref{saddle2} is true (by applying optimality condition) $\exists \; g_\gamma^\ast \in \partial^+_\gamma L(x^\ast,\gamma^\ast)$ such that 
       \begin{equation*}
       \begin{split}
          & \inprod{-g_\gamma^\ast}{\gamma-\gamma^\ast} \geq 0
          \\ <=> \;  & 0 \in  \nabla^2 R_\gamma(\gamma)^{-1} (\partial^+_\gamma L(x^\ast,\gamma^\ast) - \mathcal{N}_\mathcal{Y}(\gamma)).
           \end{split}
       \end{equation*}
       The converse can be proved in exactly the same way. 
   \end{proof}
    In the next theorem, we analyze the stability of \eqref{primal dynamics}-\eqref{dual dynamics} by exploiting the results of the Theorem \ref{Cortes}.
   \begin{theorem}
       Every equilibrium point $(x^\ast, \gamma^\ast) \in \mathcal{S}$ of \eqref{primal dynamics}-\eqref{dual dynamics} is Lyapunov stable. Furthermore,  every Caratheodory solution of the \eqref{primal dynamics}-\eqref{dual dynamics} with  initial condition $(x_0,\gamma_0) \in \mathcal{X} \times \mathcal{Y}$ asymptotically converges to the set  $\mathcal{S}$.
       \label{sability}
\end{theorem}
\begin{proof}
    Choose any $(x^\ast,\gamma^\ast) \in \mathcal{S}$. Then define the Lyapunov function $V: \mathcal{X} \times \mathcal{Y} \to [0,\infty)$ as
    \begin{equation*}
        V(x,\gamma) = \mathbb{D}_{R_x}(x^\ast,x) + \mathbb{D}_{R_\gamma}(\gamma^\ast,\gamma).
    \end{equation*}
    The term $\mathbb{D}_{R_x}(x^\ast,x)$ and $\mathbb{D}_{R_\gamma}(\gamma^\ast,\gamma)$ are referred to  as Bregman Divergence associated with the strongly convex functions $R_x$ and $R_\gamma$ and are defined as follows:
    \begin{equation*}
        \mathbb{D}_{R_x}(x^\ast,x) = R_x(x^\ast) - R_x(x) - \inprod{\nabla R_x(x)}{x^\ast-x}.
    \end{equation*}
    and
    \begin{equation*}
        \mathbb{D}_{R_\gamma}(\gamma^\ast,\gamma) = R_\gamma(\gamma^\ast)  - R_\gamma(\gamma) - \inprod{\nabla R_\gamma(\gamma)}{\gamma^\ast-\gamma}.
    \end{equation*}
    Because of $\sigma_R$-strong convexity of $R_x$ and $R_\gamma$, we have
    \begin{equation*}
        V(x,\gamma) \geq \sigma_R (\norm{x-x^\ast}^2+ \norm{\gamma-\gamma^\ast}^2)   \; \; \text{and} \; \; V(x^\ast,\gamma^\ast) = 0.
    \end{equation*}
    The set valued Lie derivative of $V$ with respect to \eqref{primal dynamics}-\eqref{dual dynamics} are as follows:
    \begin{equation}
        \begin{split}
            \Tilde{\mathcal{L}} V( x,\gamma) = \{ \inprod{\nabla V(x,\gamma)}{(\zeta_1,\zeta_2)} \; | &  \; \zeta_1 \in -\nabla^2 R_x(x)^{-1}(\partial_x L(x,\gamma) + \mathcal{N}_\mathcal{X}(x) ) \; ,
             \; \zeta_2 \in \nabla^2 R_\gamma(\gamma)^{-1} (\partial^+_\gamma L(x,\gamma) - \mathcal{N}_\mathcal{Y}(\gamma)) \} .
        \end{split}
        \label{Lie}
    \end{equation}
    Consider any $\zeta_1 \in -\nabla^2 R_x(x)^{-1}(\partial_x L(x,\gamma) + \mathcal{N}_\mathcal{X}(x)$ and $ \zeta_2 \in \nabla^2 R_\gamma(\gamma)^{-1} (\partial^{+}_\gamma L(x,\gamma) - \mathcal{N}_\mathcal{Y}(\gamma))$ then $\exists \; 
 g_x \in \partial_x L(x,\gamma)$ and $\eta_x \in \mathcal{N}_\mathcal{X}(x)$ such that
 \begin{equation*}
     \zeta_1 = - \nabla^2R_x(x)^{-1} (g_x+\eta_x)
 \end{equation*}
 and similarly,
  \begin{equation*}
     \zeta_2 =  \nabla^2R_\gamma(\gamma)^{-1} (g_\gamma-\eta_\gamma).
 \end{equation*}
 On the other hand from the definition of $V$
 \begin{equation*}
     \begin{split}
         & \nabla_x V(x,\gamma) = \nabla^2 R_x(x) (x-x^\ast)
         \\ & \nabla_\gamma V(x,\gamma) = \nabla^2 R_\gamma(\gamma) (\gamma-\gamma^\ast).
     \end{split}
 \end{equation*}
 Hence, we get 
 \begin{equation}
     \begin{split}
          \inprod{\nabla V(x,\gamma)}{(\zeta_1,\zeta_2)} 
          = & \inprod{-g_x-\eta_x}{x-x^\ast} + \inprod{\gamma-\gamma^\ast}{g_\gamma-\eta_\gamma} 
         \\ \leq & L(x^\ast,\gamma) - L(x,\gamma) + L(x,\gamma)- L(x,\gamma^\ast) \leq  0. 
     \end{split}
     \label{mimri}
 \end{equation}
The last inequality in \eqref{mimri} is a consequence of the definition of a saddle point. On the other hand, the first inequality holds from the (sub/super) gradient of the function $L$ and  also by the definition of the normal cone.
\\ Thus, $\max \Tilde{\mathcal{L}} V(x,\gamma) \leq 0$ $\forall \; (x,\gamma) \in \mathcal{X}\times \mathcal{Y}$. Hence, every $(x^\ast,\gamma^\ast) \in \mathcal{S}$ are Lyapunov stable for \eqref{primal dynamics}-\eqref{dual dynamics}. Next, we have to prove 
that the set $\mathcal{S}$ is asymptotically stable. In other words, all Caratheodory solution of \eqref{primal dynamics}-\eqref{dual dynamics} for any initial condition $(x_0,\gamma_0) \in \mathcal{X}\times \mathcal{Y}$ asymptotically converge to  the set $\mathcal{S}$. For that consider the same Lyapunov function $V$ as already assumed .

Assume that $V(x_0,\gamma_0) = c$  $> 0$. Then from the preceding discussions, we get that all Caratheodory solution of the ordinary differential inclusions \eqref{primal dynamics}-\eqref{dual dynamics} must lie within the set $\{ (x,\gamma) | 
 \; V(x,\gamma) \leq c\} \cap \mathcal{X} \times \mathcal{Y} (:=M_c)$.  This set is compact and also strongly invariant. 
 \\ Hence, by applying La-Salle's invariance principle all the Caratheodory's solutions of \eqref{primal dynamics}-\eqref{dual dynamics} must converge to largest weakly invarint set $M$ contained in 
  $M_c \cap \{(x,\gamma)\in \mathcal{X} \times \mathcal{Y} \; | \; 0 \in \Tilde{\mathcal{L}}V(x,\gamma) \}$. However, $0 \in \Tilde{\mathcal{L}}V(x,\gamma) \}$ imply from \eqref{Lie} and  \eqref{mimri} that $\exists \; (g_x,g_\gamma) \in \partial_x L(x,\gamma) \times \partial^+_\gamma L(x,\gamma)$ and $(\eta_x,\eta_\gamma) \in \mathcal{N}_\mathcal{X}(x)\times \mathcal{N}_\mathcal{Y}(\gamma)$ such that 
 \begin{equation*}
     0 = \inprod{x-x^\ast}{-g_x-\eta_x} + \inprod{\gamma-\gamma^\ast}{g_\gamma-\eta_\gamma}.
 \end{equation*}
 On applying the definition of normal cone 
 \begin{equation}
     \inprod{g_x}{x-x^\ast}-\inprod{g_\gamma}{\gamma-\gamma^\ast} \leq 0.
     \label{31}
 \end{equation}
 Claim - Inequality \eqref{31} directly implies $(x,\gamma) \in \mathcal{S}$.
 \\ This could be shown in the following way. By the definition of subgradient
 \begin{equation}
     \inprod{g_x}{x-x^\ast} \geq L(x,\gamma) - L(x^\ast,\gamma)
     \label{32}
 \end{equation}
 and from  supergradient 
 \begin{equation}
    - \inprod{g_\gamma}{\gamma-\gamma^\ast}  \geq -L(x,\gamma) + L(x,\gamma^\ast).
     \label{33}
 \end{equation}
 Adding both \eqref{32} and \eqref{33} and from \eqref{31} we get 
 \begin{equation*}
     0 \geq \inprod{g_x}{x-x^\ast}   - \inprod{g_\gamma}{\gamma-\gamma^\ast} \geq L(x,\gamma^\ast) - L(x^\ast,\gamma) \geq 0. 
 \end{equation*}
 The final inequality arises directly from the definition of saddle points. Such a scenario occurs only if $L(x,\gamma^\ast) = L(x^\ast,\gamma)$. Considering Assumption \ref{saddle assumption}, we can conclude that $(x,\gamma)$ $\in$ $\mathcal{S}$.
 
 Hence every equilibrium point $(x^\ast,\gamma^\ast) \in \mathcal{S}$ is  Lyapunov stable and every Caratheodory solution of \eqref{primal dynamics} and \eqref{dual dynamics} starting from any initial condition asymptotically converges to $\mathcal{S}$. 
\end{proof}

In this section, we explored the continuous-time dynamics of saddle point mirror descent, examining them through the perspectives of both the tangent cone and the normal cone. A key consideration emerges: how do we effectively approximate these dynamics in real-world scenarios? While conventional numerical integration schemes, as outlined in \cite{projected}, offer a viable solution, our work takes a step further. In the subsequent section, we present a more generalized version – a SSPMD algorithm – designed to address the problem stated in \eqref{primal}. Additionally, we establish a connection between the SSPMD algorithm and the projected dynamical system in the non-Euclidean domain discussed thus far.  It is worth emphasizing that we leverage the analytical insights from this section to elucidate the convergence analysis of the SSPMD algorithm.
\section{SSPMD Algorithm}
To initiate this section, we make the foundational assumption of having an oracle that provides a stochastic (sub/super) gradient denoted as $G(x,\gamma,\zeta)$ for the function $L(x,\gamma)$ at a specific point $(x,\gamma)$. This stochastic (sub/super) gradient is structured as follows:
\begin{equation*}
    G(x,\gamma,\zeta) =  \begin{bmatrix}
        G_x(x,\gamma,\zeta) \\
        G_\gamma(x,\gamma,\zeta). 
    \end{bmatrix}
\end{equation*}
Here, $\zeta$ represents a random vector defined on a probability space $(\Omega, \mathcal{F}, \mathbb{P})$.  This scenario finds its motivation when the objective function in \eqref{primal} takes the form $\mathbb{E}[F(x,\gamma,\zeta)]$. In such instances, obtaining the exact (sub/super) gradient of the function becomes computationally challenging, particularly for high-dimensional problems. In some cases, it might even be infeasible, especially when the probability distribution of $\zeta$ is unknown. Consequently, we rely on the stochastic (sub/super) gradient of the function $L(x,\gamma)$ to navigate these challenges. 

In this case, for each $(x,\gamma) \in \mathcal{X} \times \mathcal{Y}$, we initially select a sample point $\zeta(\omega)$, where $\omega \in \Omega$, and subsequently compute the (sub/super) gradient $G(x,\gamma,\zeta)$ of the function $F(x,\gamma,\zeta(\omega))$. Thus, following the definition of the (sub/super) gradient:

\begin{equation}
F(x_1,\gamma,\zeta(\omega)) \geq F(x,\gamma,\zeta(\omega)) + \langle G_x(x,\gamma,\zeta(\omega)), x_1 - x \rangle \quad \forall \; x_1 \in \mathcal{X}
\label{subg1}
\end{equation}
and
\begin{equation}
- F(x,\gamma_1,\zeta(\omega)) \geq - F(x,\gamma, \zeta(\omega)) - \langle G_\gamma (x,\gamma,\zeta(\omega)),\gamma_1 - \gamma \rangle \quad \forall \; \gamma_1 \in \mathcal{Y}.
\label{superg}
\end{equation}
By taking the expectation on both sides of \eqref{subg1} and \eqref{superg}, we can establish that:
\begin{equation*}
\mathbb{E}\Big{[} \begin{bmatrix}
        G_x(x,\gamma,\zeta)
        \\ G_\gamma(x,\gamma,\zeta) 
    \end{bmatrix} \Big{]}  =  \begin{bmatrix}
        g_x(x,\gamma) \in \partial_x L(x,\gamma)
        \\ g_\gamma(x,\gamma) \in \partial^+_\gamma L(x,\gamma) 
    \end{bmatrix}.
\end{equation*}
Utilizing these stochastic (sub/super) gradients, we  present the SSPMD algorithm. Prior to that, let's define two functions $\Bar{R}_x : \mathbb{R}^n \to \Bar{R}$ and $\Bar{R}_\gamma : \mathbb{R}^m \to \Bar{R}$ as follows:
 \begin{equation*}
    \Bar{R}_x(x) = R_x(x) + \mathcal{I}_\mathcal{X}(x) \; \; \; \text{and} \; \; \; \Bar{R}_\gamma(\gamma) = R_\gamma(\gamma) + \mathcal{I}_\mathcal{Y}(\gamma).
\end{equation*}
From this point onward, we make the assumption that both $R_x$ and $R_\gamma$ are $\rho_R$ smooth.

\textbf{SSPMD Algorithm:}
At $(n= 0)$ consider the initial iterates $x(0)$ and $\gamma(0)$. Then the update for the variable $x(n)$ is
\begin{equation}
    \begin{split}
        & y(n+1) = \nabla R_x(x(n)) - \alpha(n) G_x(x(n),\gamma(n),\zeta)
        \\ & x(n+1) = \nabla \Bar{R_x}^\ast (y(n+1))
    \end{split}
    \label{primal update}
\end{equation}
and the update for dual variable is 
\begin{equation}
    \begin{split}
        & \delta(n+1) = \nabla R_\gamma (\gamma(n)) + \alpha(n) G_\gamma(x(n),\gamma(n),\zeta)
        \\ & \gamma(n+1) = \nabla \Bar{R_\gamma}^\ast (\delta(n+1)).
    \end{split}
    \label{dual update}
\end{equation}
Here, $y(n+1)$ and $\delta(n+1)$ serve as intermediate variables, solely employed to compute the subsequent iterates $x(n+1)$ and $\gamma(n+1)$ with no additional significance. The parameter $\alpha(n) > 0$ is the step-size. Consequently, from \eqref{primal update} and \eqref{dual update}, we obtain a discrete-time stochastic process $(x(n),\gamma(n))$, and our focus is on determining whether $(x(n),\gamma(n)) \to \mathcal{S}$ almost surely.
Before presenting the almost sure convergence result, we establish some necessary assumptions. To facilitate this, we define the filtration $\mathcal{F}_n$ with respect to the stochastic process $\{x(n),\gamma(n)\}$ as follows:
\begin{equation*}
    \mathcal{F}_n = \sigma (\{x(i),\gamma(i)\}, \; \;  0 \leq i \leq n\}
\end{equation*}
\begin{assumption}
\begin{equation*}
        \mathbb{E}[G_x(n)|\mathcal{F}_n] = g_x(n) \; \; \text{a.s.}   \; \; \text{where,} \; \;  g_x(n) \in \partial_x L(x(n), \gamma(n))
    \end{equation*}
    \begin{equation*}
         \mathbb{E}[G_\gamma(n)|\mathcal{F}_n] = g_\gamma(n) \; \; \text{a.s.}  \; \; \text{where,} \; \; g_\gamma(n) \in \partial^+_\gamma L(x(n), \gamma(n)).
    \end{equation*}
    \label{expectation}
 \end{assumption}
  Henceforth, we  adopt this abbreviated notation $G_x(n)$ in place of  $G_x(x(n),\gamma(n),\zeta)$   for the sake of brevity. Similarly, we use $G_\gamma(n)$. 
  \begin{assumption} There exists $K > 0$ such that a.s.
      \begin{equation*}
         \mathbb{E}[\norm{G_x(n)}^2| \mathcal{F}_n] \leq K
     \end{equation*}
     and
      \begin{equation*}
         \mathbb{E}[\norm{G_\gamma(n)}^2| \mathcal{F}_n] \leq K.
     \end{equation*}
     \label{variance}
  \end{assumption}
  \begin{assumption} The step size $\alpha(n)$ satisfy the following Assumption.
      \begin{equation*}
          \sum\limits_{n \geq 1} \alpha(n) = \infty \; \; \text{and} \; \; \sum\limits_{n \geq 1} \alpha(n)^2 < \infty.
      \end{equation*}
      \label{steps}
  \end{assumption}
 Under these assumptions, our objective is to establish that the sequence $\{x(n),\gamma(n)\} \in \mathcal{X} \times \mathcal{Y}$ converges almost surely to $\mathcal{S}$. Before proceeding with this, we need to introduce the following technical lemma, which provides a continuous-time analog of the iterate updates in the SSPMD algorithm, as given in \eqref{primal update}-\eqref{dual update}.
 \begin{theorem}
    Let $R : \mathbb{R}^n \to \mathbb{R}$ be a strongly convex function, and $\Bar{R}(x) = R(x) + \mathcal{I}_\mathcal{X}(x)$, where $\mathcal{X} \subseteq \mathbb{R}^n$ is a closed convex set. Assuming $x \in \mathcal{X}$ and $d \in \mathbb{R}^n$, let 
    \begin{equation*}
        x_1 = \nabla \Bar{R}^\ast(\nabla R(x) + \alpha d),
    \end{equation*}
    
    then the following result holds
    \begin{equation}
        \lim_{\alpha \downarrow 0} \frac{\nabla \Bar{R}^\ast(\nabla R(x) + \alpha d) - \nabla \Bar{R}^\ast(\nabla R(x))}{\alpha} = \mathcal{P}_{\mathcal{T}_\mathcal{X}(x)}^{x} (\nabla^2 R(x)^{-1} d).
        \label{shdh}
    \end{equation}
    \label{rumsia}
\end{theorem}
\begin{proof}
    Consider the mirror descent updates 
    \begin{equation}
        y = \nabla R(x) + \alpha d
        \label{M1}
    \end{equation}
    and 
    \begin{equation}
        x_1 = \nabla \Bar{R}^\ast(y).
        \label{M2}
    \end{equation}
    From \eqref{M2}, we get that 
    \begin{equation*}
        x_1 = \argmin\limits_{z \in \mathcal{X}} \{ \inprod{z}{y} - R(z) \}.
    \end{equation*}
    From the first-order optimality  condition $\exists \; \eta_x \in \mathcal{N}_\mathcal{X}(x_1)$ such that 
    \begin{equation*}
        y - \nabla R(x_1) + \eta_x  = 0.
    \end{equation*}
    From \eqref{M1}, we get that
    \begin{equation*}
        \nabla R(x_1) = \nabla R(x) + \alpha (d + \eta_x' )
    \end{equation*}
    where, $\eta_x' = \frac{\eta_x}{\alpha} \in \mathcal{N}_\mathcal{X}(x_1)$.
   Hence, 
   \begin{equation*}
       x_1 = (\nabla R)^{-1} ( \nabla R(x) + \alpha (d + \eta_x' ) ).
   \end{equation*}
   It can be seen that $(\nabla R)^{-1}(.)$ is continuously differentiable, hence by using Taylor series expansion 
   \begin{equation}
       x_1 = (\nabla R)^{-1} (\nabla R(x)) + D (\nabla R)^{-1} (\nabla R(x)) (\alpha( d +\eta_x'))  + \small{o} (\alpha)
       \label{37}
   \end{equation}
   From Inverse Function Theorem it could be verified  that $\forall \; d \in \mathbb{R}^n$ 
   \begin{equation}
       D (\nabla R)^{-1} (\nabla R(x)) (d) = D \nabla R((\nabla R)^{-1} \nabla R(x)) ^{-1} (d) = \nabla^2  R(x)^{-1} d.
       \label{momo}
   \end{equation}
Hence, from \eqref{37} we get that
\begin{equation*}
    x_1 = x + \alpha (\nabla^2 R (x))^{-1} (d +\eta_x') + \small{o} (\alpha).
\end{equation*}
 From the optimality condition,  \eqref{momo} can be  equivalently written as  
\begin{equation}
    x_1 = \argmin\limits_{\nu \in \mathcal{X}} \norm{\nu - x - \alpha \nabla^ 2 R(x)^{-1}(d)}_x^2 + \small{o} (\alpha).
    \label{39}
\end{equation}
 Ignore $\small{o}(\alpha)$ for now  \eqref{39} is true if and only if 
\begin{equation}
    \inprod{- \nabla^2 R (x)(x_1-x -\alpha \nabla^2 R(x)^{-1} d)}{y - x_1} \leq 0 \; \; \forall \; \; y \in \mathcal{X} .
    \label{rumpi}
\end{equation}
In other words, by reorganizing \eqref{rumpi} we get
\begin{equation}
     \inprod{-\nabla^ 2 R(x) \Big{(} \frac{x_1 -x}{\alpha} -\nabla^2 R(x)^{-1} d \Big{)}} {\frac{y-x}{\alpha} - \frac{x_1-x}{\alpha}} \leq 0.
    \label{40}
\end{equation}
%The inequality in \eqref{40} implies 
%\begin{equation*}
 %   \frac{x(n+1)-x(n)}{\alpha} = \mathcal{P}_{\frac{\mathcal{X}-\{x(n)\}}{\alpha}}^{x(n)} (\nabla^2 R(x(n))^{-1} d) .
%\end{equation*}
If we let $\alpha \downarrow 0$ in \eqref{40}, then $\forall \; \nu \in \Lim\limits_{\alpha \downarrow 0} \frac{\mathcal{X}-\{x\}}{\alpha}  = \mathcal{C}$, we have
\begin{equation}
     \inprod{-\nabla^ 2 R(x) \Big{(} \frac{x_1 -x}{\alpha} -\nabla^2 R(x)^{-1} d \Big{)}} {\nu - \frac{x_1-x}{\alpha}} \leq 0.
     \label{bimla}
\end{equation}
Since $\mathcal{T}_\mathcal{X}(x) = \Bar{\mathcal{C}}$, $\forall \;  \nu \in \mathcal{T}_\mathcal{X}(x)$, then  \eqref{bimla} holds. Thus, from the optimality condition, we obtain 
\begin{equation*}
    \lim\limits_{\alpha \downarrow 0} \frac{x_1-x}{\alpha} = \mathcal{P}_{\mathcal{T}_\mathcal{X}(x)}^{x} (\nabla^2 R(x)^{-1} d).
\end{equation*}
Note that since $\lim\limits_{\alpha \downarrow 0} \frac{\small{o} (\alpha)}{\alpha} = 0$, hence \eqref{shdh} holds.
\end{proof}
 \begin{remark}
   Theorem \ref{rumsia} provides the directional derivative of $\nabla \Bar{R}^\ast (z)$, which is commonly referred to in the literature as stating that $\nabla \Bar{R}^\ast (z)$ is Gateaux semidifferentiable \cite{10.5555/3384669} for all $z \in \nabla R(\mathcal{X})$ . The limit in Theorem \ref{rumsia} is denoted as $d \nabla \Bar{R}^\ast(z;d)$ and is  equal to the projection of the vector $d$ onto the tangent cone in the non-Euclidean metric.
\end{remark}
 \begin{remark}
    Observe the expression \eqref{39} within the proof of Theorem \ref{rumsia}, which highlights a crucial insight. The algorithmic steps of the mirror descent algorithm can be interpreted as an approximate application of the projected subgradient method within the context of a Riemannian manifold. This observation is particularly noteworthy, as it generalizes a similar result previously established  for unconstrained optimization problem \cite{7004065}. However, the significance of this result lies in its broader applicability, encompassing a more general class of constrained optimization problem. 
\end{remark}

The importance of Theorem \ref{rumsia} suggests that the proposed projected dynamical system in a non-Euclidean domain can be approximated numerically using the steps outlined in the saddle point mirror descent algorithm, as given in equations \eqref{primal update}-\eqref{dual update} (with $G_x(n) = g_x(n) \in \partial_x L(x(n),\gamma(n))$ and $G_\gamma(n) = g_\gamma(n)  \in \partial^+_\gamma L(x(n),\gamma(n))$). As we proceed through the paper, it becomes evident that we delve into the almost sure convergence analysis of a more generalized version - the SSPMD algorithm, and in the subsequent section, the SZSPMD algorithm. Here, Theorem \ref{rumsia} serves as a link between the iterative steps of the SSPMD algorithm, outlined in equations \eqref{primal update}-\eqref{dual update}, and the projected dynamical systems presented in equations \eqref{projected primal}-\eqref{projected dual}. This connection is further clarified in  the following Corollary.

\begin{Corollary}
    The iterates of the stochastic mirror descent algorithm as given in  \eqref{primal update}-\eqref{dual update} can equivalently be   written in the following way:
    \begin{equation}
    \begin{split}
        x(n+1) & =  x(n) +  \alpha(n)\mathcal{P}_{\mathcal{T}_\mathcal{X}(x(n))}^{x(n)} (\nabla^2 R_x(x(n))^{-1} (-G_x(n))) + \small{o} (\alpha (n))
        \\ & = x(n) + \alpha(n) \nabla^ 2 R_{x}(x(n))^{-1} (-G_x(n) - \eta_x(n)))  + \small{o}(\alpha(n)).
       \end{split}
       \label{43}
    \end{equation}
    and similarly for dual variable $\gamma(n)$ we have
     \begin{equation}
    \begin{split}
        \gamma(n+1) & =  \gamma(n) +  \alpha(n) \mathcal{P}_{\mathcal{T}_\mathcal{Y}(\gamma(n))}^{\gamma(n)} (\nabla^2 R_\gamma(\gamma(n))^{-1} (G_\gamma(n))) + \small{o} (\alpha (n))
        \\ & = \gamma(n) + \alpha(n) \nabla^ 2 R_{\gamma}(\gamma(n))^{-1} (G_\gamma(n) - \eta_\gamma(n))  + \small{o}(\alpha(n)).
       \end{split}
       \label{44}
    \end{equation}
    \label{cor}
\end{Corollary} 
Before delving further, it's worth noting that the stochastic (sub/super) gradient $G(n)$ can also be expressed in the following manner:
\begin{equation*}
    G_x(n) = g_x(n)  + M^1(n+1).
\end{equation*}
Here, $g_x(n) \in \partial_x L(x(n), \gamma(n))$, and $\{M^1(n)\}$ forms a martingale difference sequence with respect to the filtration ${\mathcal{F}_n }$. Similarly,
\begin{equation*}
    G_\gamma(n) = g_\gamma(n) + M^2(n+1).
\end{equation*}
In this context, equations \eqref{43}-\eqref{44} can be rewritten in the following way
\begin{equation*}
    x(n+1) = x(n) + \alpha(n) (\nabla^2 R_x(x(n))^{-1}(-g_x(n)  - \eta_x(n))) + M^1_x (n+1) + \small{o}(\alpha(n))
\end{equation*}
 and 
 \begin{equation*}
     \gamma(n+1) = \gamma(n) + \alpha(n)  (\nabla^2 R_\gamma (\gamma(n))^{-1}(g_\gamma(n) - \eta_\gamma(n))) + M_\gamma^2(n+1)  + \small{o}(\alpha(n)).
 \end{equation*}
 Here, $M_x^1(n+1) = \nabla^2 R_x(x(n))^{-1} M^1(n+1)$  and $M_\gamma^2 (n+1) = \nabla^2 R_\gamma(\gamma(n))^{-1} M^2 (n+1)$. It is straightforward to check that both $\{M_x^1(n)\}$ and $\{M_\gamma^2(n)\}$ form a martingale difference sequence with respect to the filtration $\mathcal{F}_n$. 

  Next, we define the continuous-time interpolated trajectories $({\Hat{x}(t),\Hat{\gamma}(t)})$ corresponding to the iterates of the stochastic mirror descent algorithm $({x(n),\gamma(n)})$ in \eqref{primal update}-\eqref{dual update} as follows:
\begin{equation}
    \Hat{x}(t) = x(n) + (x(n+1)-x(n)) \frac{t-t(n)}{t(n+1)-t(n)} \; \; \; \forall \; \; t \in I_n
\label{pinter}
\end{equation}
where, $I_n = [t(n),t(n+1)]$.
In a similar fashion, we can define the interpolated trajectory for $\gamma(t)$ as follows
\begin{equation}
    \Hat{\gamma}(t) = \gamma(n) + (\gamma(n+1)-\gamma(n)) \frac{t - t(n)}{t(n+1)-t(n)} \; \; \; \forall \; \; t \in I_n.
    \label{dinter}
\end{equation}
The preceding discussion in this section culminates in the following Proposition.
\begin{proposition}
    The interpolated trajectory $( \Hat{x}(t), \Hat{\gamma}(t) )$  is  an asymptotic pseudo trajectory of the differential inclusions \eqref{primal dynamics}-\eqref{dual dynamics} (or equivalently, asymptotic pseudo trajectory of \eqref{projected primal}-\eqref{projected dual}). In other words for any $T >0 $ (a.s.)
    \begin{equation*}
        \lim\limits_{t \to \infty} \;  \inf\limits_{x(s) \in S_{\hat{x}(t)}} \; \sup\limits_{0 \leq s \leq T} \norm{\Hat{x}(t+s)-x(s)} = 0  \; \; \text{and}
    \end{equation*}
    \begin{equation*}
        \lim\limits_{t \to \infty} \;  \inf\limits_{\gamma(s) \in S_{\hat{\gamma}(t)}} \; \sup\limits_{0 \leq s \leq T} \norm{\Hat{\gamma}(t+s)-\gamma(s)} = 0  
    \end{equation*}
    where, $S_{\Hat{x}(t)}$ and  $S_{\Hat{\gamma}(t)}$ are  the set of all solutions $(x(s),\gamma(s))$ of the differential inclusions \eqref{primal dynamics}-\eqref{dual dynamics} with initial conditions $x(0) = \Hat{x}(t)$ and $\gamma(0) = \Hat{\gamma}(t)$. 
    \label{benaim1}
\end{proposition}
By leveraging Proposition \ref{benaim1}, we establish the almost sure convergence of the iterates of the SSPMD algorithm.
\subsection{Convergence Analysis of iterates of SSPMD Algorithm with unbiased stochastic (sub/super) gradient}
In this subsection, we establish the almost sure convergence of the iterates of the SSPMD algorithm under Assumption \ref{expectation}. Furthermore, in the subsequent section, we extend the convergence analysis beyond Assumption \ref{expectation}, particularly in the context of the SZSPMD algorithm.
\begin{theorem}
    The iterate sequence generated by SSPMD algorithm \eqref{primal update}-\eqref{dual update}  converge to $\mathcal{S}$ almost surely. 
    \label{almstsure}
\end{theorem}
\begin{proof}
We aim to demonstrate the almost sure convergence of $(\Hat{x}(t),\Hat{\gamma}(t) )$ to $\mathcal{S}$, which, in turn, implies the almost sure convergence of $(x(n),\gamma(n))$ to $\mathcal{S}$. In other word it means that  $\forall \; \epsilon >0 $, $\exists \; T_0 > 0 $ such that $\forall \; t \geq T_0$ we have $(\Hat{x}(t), \Hat{\gamma}(t)) \in N_\epsilon(\mathcal{S}) $. 

Define $V(x,\gamma) : \mathbb{R}^n \times \mathbb{R}^m \to \mathbb{R}^{+}$ for each $(x^\ast,\gamma^\ast) \in \mathcal{S}$ (as defined in Theorem \ref{sability}) as follows:
\begin{equation*}
    V(x,\gamma) = \mathbb{D}_{R_x}(x^\ast,x) +   \mathbb{D}_{R_\gamma}(\gamma^\ast,\gamma).
\end{equation*}

However, from the definition of strong convexity of Bregman divergence, we have

\begin{equation}
 V(\Hat{x}(t),\Hat{\gamma}(t)) \geq \frac{\sigma_R}{2}(\norm{\Hat{x}(t)-x^\ast}^2 + \norm{\Hat{\gamma}(t)-\gamma^\ast}^2).
    \label{rax}
\end{equation}

From \eqref{rax}, it's evident that to demonstrate almost sure convergence, it suffices to show that  $\forall \; \epsilon > 0$, $\exists \; T_0 > 0$ such that $\forall \; t \geq T_0$ $\exists \; (x^\ast,\gamma^\ast) \in \mathcal{S}$ we have $V(\Hat{x}(t),\Hat{\gamma}(t)) \leq \epsilon$ a.s.

   We proceed with the proof in two steps. For the first step, let's assume that $(\hat{x}(t),\hat{\gamma}(t)) \in A \subseteq \mathbb{R}^n \times \mathbb{R}^m$. The set $A$ is assumed to be compact since the sequence $\{x(n),\gamma(n)\}$ is bounded almost surely.
   
   In this step, we aim to show that $\forall \; \epsilon > 0$, if $(x(t),\gamma(t))$ be any viable Caratheodory solution  of the ordinary differential inclusion \eqref{projected primal}-\eqref{projected dual} with any initial condition $(x_0,\gamma_0) \in A$, then $\exists \; T > 0$ (irrespective of the initial condition as long as they belong to $A$) and $(x^\ast,\gamma^\ast) \in \mathcal{S}$ such that $V(x(T) , \gamma(T)) \leq \frac{\epsilon}{2}$.

Define the function $V^\ast (x,\gamma) : \mathbb{R}^n \times \mathbb{R}^m \to \mathbb{R}^+$ as follows:
\begin{equation}
    V^\ast(x,\gamma) = \min\limits_{(x^\ast,\gamma^\ast) \in \mathcal{S}} \{ V(x,\gamma) =  \mathbb{D}_{R_x} (x^\ast,x)  + \mathbb{D}_{R_\gamma} (\gamma^\ast,\gamma) \}.
    \label{rama}
\end{equation}

Note that, based on the definition of Bregman Divergence, $V(x,\gamma)$ is continuously differentiable with respect to  $(x,\gamma) \in \mathcal{X} \times \mathcal{Y}$ $\forall \; (x^\ast,\gamma^\ast) \in \mathcal{S}$. Since the set $\mathcal{S}$ is compact, according to Danskin's Theorem, $V^\ast(x,\gamma)$ is also continuously differentiable.

As the set $\mathcal{S}$ is compact, there exists a minimum of \eqref{rama} for each $(x,\gamma) \in \mathbb{R}^n \times \mathbb{R}^m$, denoted $(x^\ast_1,\gamma^\ast_1) \in \mathcal{S}$. In other words:
\begin{equation*}
    V^\ast(x,\gamma) = \mathbb{D}_{R_x}(x^\ast_1,x) + \mathbb{D}_{R_\gamma}(\gamma_1^\ast,\gamma).
\end{equation*}

According to Danskin's Theorem, we have:
\begin{equation}
    \nabla V^\ast (x,\gamma) = \nabla_x \mathbb{D}_{R_x}(x^\ast_1,x) + \nabla_\gamma \mathbb{D}_{R_\gamma}(\gamma^\ast_1,\gamma)
    \label{ydu}.
\end{equation}

 Let $\alpha : = \max\limits_{(x,\gamma) \in A} V^\ast(x,\gamma)$. Consider the set $W = \{ (x,\gamma) \in \mathcal{X} \times \mathcal{Y}  | V^\ast(x,\gamma) \leq \alpha \} \subseteq A$. Clearly, the set $W$ is closed, as per \eqref{rax}.

Define a neighborhood around the set $\mathcal{S}$ in the space $\mathcal{X} \times \mathcal{Y}$ as $U_{\frac{\epsilon}{3}}(\mathcal{S}) = \{(x,\gamma) \in \mathcal{X} \times \mathcal{Y} \mid V^\ast(x,\gamma) < \frac{\epsilon}{3}\}$.

From the proof of Theorem \ref{sability}, it's evident that the set-valued Lie derivatives of $V(x,\gamma)$, denoted as $(\tilde{\mathcal{L}}V(x,\gamma))$ for any $(x^\ast,\gamma^\ast) \in \mathcal{S}$ satisfies:
\begin{equation*}
\Tilde{\mathcal{L}}V(x,\gamma) \leq L(x^\ast,\gamma)- L(x,\gamma^\ast)  \leq 0 \; \; \forall \; (x,\gamma) \in \mathcal{X} \times \mathcal{Y} .
\end{equation*}

According to \eqref{ydu} and the definition of Lie derivative, $\Tilde{\mathcal{L}} V^\ast(x,\gamma) \leq 0$. 
In other words, $\dot{V}^\ast$ along the trajectories of the ordinary differential inclusion given in \eqref{projected primal}-\eqref{projected dual} is non-positive. Therefore, $U_\frac{\epsilon}{3}$ as well as $W$ are positively invariant. 

Consider the set $W \setminus U_{\frac{\epsilon}{3}} (\mathcal{S}) \subseteq \mathcal{X} \times \mathcal{Y}$, which is compact. Thus, $\exists $ $\beta > 0$ and $\forall \; (x^\ast,\gamma^\ast) \in \mathcal{S}$ such that $L(x,\gamma^\ast) - L(x^\ast,\gamma) \leq - \beta$ $\forall \; (x,\gamma) \in W \setminus U_\frac{\epsilon}{3}$. Clearly, in view of \eqref{ydu} we obtain: 
\begin{equation}
    \max\limits_{(x,\gamma) \in A \setminus U_\frac{\epsilon}{3}} \Tilde{\mathcal{L}} V^\ast (x,\gamma) \leq - \beta
    \label{ddd}
\end{equation}

Now, let $(x(t),\gamma(t))$ be any solution of the ordinary differential inclusion given in \eqref{projected primal}-\eqref{projected dual} with initial condition $(x_0,\gamma_0) \in A$. Then $(x(t),\gamma(t)) \in W$. This implies that $V^\ast(x(0),\gamma(0)) \leq \alpha $ and $\Dot{V}^\ast(t) \leq - \beta$ along the trajectory $(x(t),\gamma(t))$ if $(x(t),\gamma(t)) \notin U_\frac{\epsilon}{3}$.   Now, it can be verified that if we choose $T > \frac{3 \alpha - \epsilon}{\beta}$, then $V^\ast(x(T),\gamma(T)) \leq \frac{\epsilon}{3}$.

With this choice of $T$, we proceed to the second step of the proof where we show that $(\hat{x}(t),\hat{\gamma}(t))$ converge to $\mathcal{S}$ almost surely.

Since $(\Hat{x}(t),\Hat{\gamma}(t))$ is the asymptotic pseudo-trajectory of the differential inclusions \eqref{primal dynamics}-\eqref{dual dynamics}, we have (a.s.):
\begin{equation*}
       \lim\limits_{t \to \infty} \inf\limits_{x(s)  \in S_{\hat{x}(t)}} \sup\limits_{0 \leq s \leq T} \norm{\Hat{x}(t+s)-x(s)} = 0 
   \end{equation*}

For any $\delta > 0$,  $\exists \; t_0(T,\delta) > 0$ such that  $\forall \; t \geq t_0$, we have (a.s.):

   \begin{equation*}
       \inf\limits_{x(s) \in S_{\hat{x}(t)}} \sup\limits_{0 \leq s \leq T} \norm{\Hat{x}(t+s)-x(s)} \leq \frac{\delta}{2} .
   \end{equation*}
 
   Fix any $t \geq t_0$, there exists a solution $(x_t(s),\gamma_t(s))$ of the differential inclusions  \eqref{primal dynamics}-\eqref{dual dynamics}  with initial condition $x_t(0) = \Hat{x}(t)$ and $\gamma_t(0) = \Hat{\gamma}(t)$ such that
   \begin{equation}
       \sup\limits_{0 \leq s \leq T} \norm{\Hat{x}(t+s)-x_t(s)} \leq \delta \; \; \text{a.s.}
       \label{24}
   \end{equation} 
 \begin{equation}
       \sup\limits_{0 \leq s \leq T} \norm{\hat{\gamma}(t+s)- \gamma_t(s)} \leq \delta \; \; \text{a.s}.
       \label{251}
   \end{equation}

From the definition of Bregman divergence, we have: 
\begin{equation}
         \begin{split}
             & \mathbb{D}_{R_x}(x^\ast, \Hat{x}(t+T))
             = \mathbb{D}_{R_x}(x^\ast, x_t(T)) + \mathbb{D}_{R_x}(x_t(T),\Hat{x}(t+T))
            + \inprod{\nabla R_x(\Hat{x}(t+T))- \nabla R_x(x_t(T)))}{x_t(T)-x^\ast}.
         \end{split}
         \label{mimi}
     \end{equation}

 Notice that  $\mathbb{D}_{R_x}(x_t(T),\Hat{x}(t+T)) \leq \frac{\rho_R}{2} \norm{x_t(T)-\Hat{x}(t+T)}^2$ due to the $\rho_R$-smoothness of the function of $R_x$. This smoothness also provides
\begin{equation*}
\begin{split}
    \inprod{\nabla R_x(\Hat{x}(t+T))- \nabla R_x(x_t(T)))}{x_t(T)-x^\ast}
     \leq & \norm{\nabla R_x(\Hat{x}(t+T))- \nabla R_x(x_t(T)))}_\ast  \norm{{x_t(T)-x^\ast}}
    \\ \leq & \rho_R \norm{\Hat{x}(t+T) - x_t(T)}_\ast \norm{x_t(T)-x^\ast}.
\end{split}    
\end{equation*}

From \eqref{251}, we conclude that $\lVert \Hat{x}(t+T)- x_t(T) \rVert \leq \delta$ a.s. Thus, combining these results in \eqref{mimi}, for all $t \geq t_0$, we obtain (a.s.):
     \begin{equation}
         \begin{split}
 \mathbb{D}_{R_x}(x^\ast, \Hat{x}(t+T))
              \leq  \mathbb{D}_{R_x}(x^\ast, x_t(T)) + \frac{\rho_R}{2} \delta^2 + \rho_R \delta \norm{x_t(T) - x^\ast}.
         \end{split}
    \label{inequality involving Bergman}
     \end{equation}

Similarly, we can derive an analogous inequality such that  $\exists \; t_1$ such that $\forall \; t \geq t_1$
  \begin{equation}
         \begin{split}
 \mathbb{D}_{R_\gamma}(\gamma^\ast, \Hat{\gamma}(t+T))
              \leq  \mathbb{D}_{R_\gamma}(\gamma^\ast, \gamma_t(T)) + \frac{\rho_R}{2} \delta^2 + \rho_R \delta \norm{\gamma_t(T) - \gamma^\ast}.
         \end{split}
    \label{inequality involving Bergmanr}
     \end{equation}

Choosing $t_2 = \max \{t_0,t_1\}$ and adding inequalities \eqref{inequality involving Bergman} and \eqref{inequality involving Bergmanr}, we get $\forall \; t \geq t_2$ 
\begin{equation}
         \begin{split}
              V(\Hat{x}(t+T), \Hat{\gamma}
             (t+T))
              \leq   V(x_t(T),\gamma_t(T))  + \rho_R \delta^2 +  \rho_R \delta (\norm{x_t(T)-x^\ast} + \norm{\gamma_t(T) - \gamma^\ast}) \; \; \text{a.s}. 
         \end{split}
         \label{22}
     \end{equation}

   Based on Lyapunov's definition of stability and asymptotic set stability, since $(x_t(T),\gamma_t(T))$ is a solution of \eqref{projected primal}-\eqref{projected dual} with initial condition $(x_t(0),\gamma_t(0)) \in A$, there exists a constant $C > 0$ such that $\lVert x_t(T)-x^\ast \rVert + \lVert \gamma_t(T) - \gamma^\ast \rVert \leq C$ for all $t \geq 0$ and for all $(x^\ast,\gamma^\ast) \in \mathcal{S}$. Thus, according to \eqref{22}, we obtain $\forall \; (x^\ast,\gamma^\ast) \in \mathcal{S}$ (a.s.):
\begin{equation}
    V(\Hat{x}(t+T),\Hat{\gamma}(t+T)) \leq V(x_t(T),\gamma_t(T)) + \rho_R \delta ^2 + \rho_R \delta C.
    \label{abj}
\end{equation}

Since $\delta$ is free, $\rho_R \delta^2 + l \delta C \leq \frac{\epsilon}{2}$. As shown in step 1, it has been demonstrated that $V^\ast (x_t(T),\gamma_t(T)) \leq \frac{\epsilon}{3}$. Therefore, from \eqref{abj} we conclude that $\exists \; t_2 > 0$ such that $\forall \; t \geq t_2 + T$, $\exists \; (x^\ast,\gamma^\ast) \in \mathcal{S}$ where (a.s.) $V(\Hat{x}(t+T), \Hat{\gamma}(t+T)) \leq \epsilon$. 
\end{proof}
In this section, we present the almost sure convergence analysis of the iterates of the SSPMD algorithm for unbiased (sub/super) subgradients. However, in the next subsection, we relax this unbiased assumption and conduct the convergence analysis using biased (sub/super) gradients.
\subsection{Convergence Analysis  of SZSPMD Algorithm}
This subsection commences with the premise of the existence of an oracle providing the function value $L(x,\gamma)$ at any desired point $(x,\gamma) \in \mathcal{X} \times \mathcal{Y}$. This premise arises from the practical reality that in many applications, there may not be any closed-form expression available for  $L(x,\gamma)$, or even if available, computing the gradient of the function $L(x,\gamma)$ could be exceedingly computationally challenging due to large dimensions. In this context, we introduce the SZSPMD algorithm.

Within this subsection, $L(x,\gamma)$ is presumed to be continuously differentiable, leading to a reliance on approximated gradients of the function $L(x,\gamma)$ for presenting the SSPMD algorithm. We employ Nesterov's Gaussian approximation methods \cite{Nesterov2015RandomGM} to approximate the gradient of the function $L(x,\gamma)$ at any point $(x,\gamma) \in \mathcal{X} \times \mathcal{Y}$ . The algorithmic steps of the SZSPMD algorithm are as follows:

At iteration $n$, consider updates for primal variable and dual variables as $(x(n),\gamma(n))$ respectively. Then, to approximate the gradient of the function $L(x,\gamma)$ at the iterates $(x(n),\gamma(n))$, we first consider a normal random vector $u \in \mathbb{R}^{n \times m}$ with zero mean and identity standard deviation. The approximated
gradient ($\Tilde{\nabla} L(x,\gamma)$) at $(x(n),\gamma(n))$ can be calculated as follows:

\begin{equation}
    \Tilde{\nabla} L(x(n),\gamma(n)) = \frac{L((x(n),\gamma(n))+ \mu(n) u) - L(x(n),\gamma(n))}{\mu(n)}.
    \label{approx}
\end{equation}
Here, $\mu(n) > 0$ is a constant parameter. Using this approximated gradient, we compute the next iterates $(x(n+1),\gamma(n+1))$ as follows:
\begin{equation}
    \begin{split}
        & y(n+1) = \nabla R_x(x(n)) - \alpha(n) \Tilde{\nabla}_x L(x(n),\gamma(n))
        \\ & x(n+1) = \nabla \Bar{R_x}^\ast (y(n+1))
    \end{split}
    \label{primal updateb1}
\end{equation}
and for dual variable  
\begin{equation}
    \begin{split}
        & \delta(n+1) = \nabla R_\gamma (\gamma(n)) + \alpha(n) \Tilde{\nabla}_\gamma L (x(n),\gamma(n))
        \\ & \gamma(n+1) = \nabla \Bar{R_\gamma}^\ast (\delta(n+1))
    \end{split}
    \label{dual updateb1}
\end{equation}
In this subsection, our focus is on the almost sure convergence analysis of the SZSPMD algorithm for the updates provided in \eqref{primal updateb1} and \eqref{dual updateb1}. Before delving into the convergence analysis, it is essential to outline the properties of the approximated gradient of the function $L(x,\gamma)$, as given in \eqref{approx}. For that define the filtration $\{\mathcal{F}_n \}$ as follows:
\begin{equation*}
    \mathcal{F}_n = \sigma ( \{x(i),\gamma(i) \}, \; \;   0 \leq i \leq n )
\end{equation*}
The following lemma provides the properties of the approximated gradient of the function $L(x,\gamma)$. A detailed proof can be referenced in \cite{Nesterov2015RandomGM}.
\begin{lemma}
    \begin{equation*}
        \begin{split}
            &
            \norm{\mathbb{E}[\Tilde{\nabla} L(x(n),\gamma(n)) | \mathcal{F}_n] - \nabla L(x(n),\gamma(n))} = \mathcal{O}(\mu(n)) \; \; \text{a.s.} 
            \\ & \mathbb{E}[\norm{\Tilde{\nabla}L (x(n),\gamma(n)}^2 | \mathcal{F}_n] \leq \mathcal{K} 
        \end{split}
    \end{equation*}
    where, $\mathcal{K} >0$ is a constant.
    \label{nesterov}
\end{lemma}

Leveraging lemma \ref{nesterov}, we can express the approximated gradient in the following manner:
\begin{equation*}
    \Tilde{\nabla}_x L(x,\gamma) = \nabla_x L(x,\gamma) + b_x(n) + M^1 (n+1)
\end{equation*}
and 
\begin{equation*}
    \Tilde{\nabla}_\gamma L(x,\gamma) = \nabla_\gamma L(x,\gamma) + b_\gamma(n) + M^2 (n+1).
\end{equation*}
Given lemma \ref{nesterov}, both $b_x(n)$ and $b_\gamma(n)$ are $\mathcal{O}(\mu(n))$.  To demonstrate the almost sure convergence of $(x(n),\gamma(n))$, we employ the same methodology as in the previous subsection. To accomplish this, we first need to identify the corresponding limiting differential inclusion.
To facilitate this analysis, we define the interpolated trajectories $\Hat{x}(t)$ and $\Hat{\gamma}(t)$ based on the equations provided in \eqref{pinter} and \eqref{dinter}, corresponding to the iterate updates $(x(n),\gamma(n))$ produced by the SZSPMD algorithm, as outlined in \eqref{primal updateb1}-\eqref{dual updateb1}.  Building upon the discussions presented thus far in this paper, particularly in the context of Theorem \ref{rumsia} and with the same reasoning provided in Corollary \ref{cor}, we can reformulate the limiting differential inclusions of \eqref{primal updateb1} and \eqref{dual updateb1} as follows:
\begin{equation}
\begin{split}
 &   \Dot{x}(t) \in -\nabla^2R_x(x(t))^{-1} (\nabla_x L(x(t),\gamma(t))  + \mathcal{N}_\mathcal{X}(x(t))) + b_x'(t)
    \end{split}
    \label{projected primalip1}
\end{equation}
\begin{equation*}
      x(0) = x_0 \in \mathcal{X}
\end{equation*}
\begin{equation}
    \Dot{\gamma}(t) \in  \nabla^2 R_\gamma(\gamma(t))^{-1}(\nabla_\gamma L(x(t),\gamma(t)) - \mathcal{N}_\mathcal{Y}(\gamma(t))) + b_\gamma'(t)  
    \label{projected dualip1}
\end{equation}
\begin{equation*}
    \gamma(0) = \gamma_0 \in \mathcal{Y}.
\end{equation*}

 Here $b_x'(t) = \nabla^2 R_x(x(t))^{-1} b_x(t)$ and $b_\gamma'(t) = \nabla^2 R_\gamma(\gamma(t))^{-1} b_\gamma(t)$ . The expressions for $b_x(t)$ and $b_\gamma(t)$ can be obtained by interpolating $b_x(n)$ and $b_\gamma(n)$ according to \eqref{pinter} and \eqref{dinter}. 
  It is evident from the discussion on Stochastic Approximation in Section \ref{mathpre}  that the interpolated trajectory $(\Hat{x}(t),\Hat{\gamma}(t))$ serves as the asymptotic pseudo-trajectory for the ordinary differential inclusions given in \eqref{projected primalip1}-\eqref{projected dualip1}. This is clarified in the following Proposition.

\begin{proposition}
    For any $T >0 $ we have (a.s.)
    \begin{equation*}
        \lim\limits_{t \to \infty}  \inf\limits_{x_t(s) \in S_{\Hat{x}(t)}}\sup\limits_{0 \leq s \leq T} \norm{\Hat{x}(t+s)-x_t(s)} = 0 
    \end{equation*}
    and
    \begin{equation*}
        \lim\limits_{t \to \infty} \inf\limits_{\gamma_t(s) \in S_{\Hat{\gamma}(t)}} \sup\limits_{0 \leq s \leq T} \norm{\Hat{\gamma}(t+s)-\gamma_t(s)} = 0  
    \end{equation*}
     where, $S_{\Hat{x}(t)}$ and  $S_{\Hat{\gamma}(t)}$ are  the set of all solutions $(x_t(s),\gamma_t(s))$ of the differential inclusions \eqref{projected primalip1}-\eqref{projected dualip1} with initial conditions $x_t(0) = \Hat{x}(t)$ and $\gamma_t(0) = \Hat{\gamma}(t)$. 
    \label{pic}
\end{proposition}

Considering $b_x(t)$ and $b_\gamma(t)$ as disturbance inputs for \eqref{projected primalip1} and \eqref{projected dualip1}, it's evident that if  $\norm{b_x(n)}$, $\norm{b_\gamma(n)} \leq \mathbf{B}$ $\forall \; n \geq 1$ then $\norm{b_x(t)}, \norm{b_\gamma(t)} \leq \mathbf{B}$.  However, by adjusting $\mu(n)$ in the algorithm, we can modify $\mathbf{B}$ according to our choice, as indicated in lemma \ref{nesterov}.
 In the subsequent theorem, utilizing Proposition \ref{pic}, we establish that by adjusting the value of $\mathbf{B}$, almost sure convergence of the iterates of the SZSPMD algorithm to any desired neighborhood of the set of saddle points is achieved.

Before presenting the main Theorem in this subsection, we require the following assumptions.

\begin{assumption}
    The functions $\nabla^2R_x(x) : \mathcal{X} \to \mathbb{R}^n \times \mathbb{R}^n$ and $\nabla^2 R_\gamma(\gamma): \mathcal{Y} \to \mathbb{R}^m \times \mathbb{R}^m$ are both $K_2$-Lipschitz.  In other words $\exists \; K_2 >0$, $\forall \; x_1,x_2 \in \mathcal{X}$ and $\gamma_1, \gamma_2 \in \mathcal{Y}$, the following holds:
    \begin{equation*}
        \norm{\nabla^2 R_x(x_1)-\nabla^2 R_x (x_2)} \leq K_2 \norm{x_1-x_2} \; \; \text{and} \; \;  \norm{\nabla^2 R_\gamma(\gamma_1)-\nabla^2 R_\gamma (\gamma_2)} \leq K_2 \norm{\gamma_1-\gamma_2}.
    \end{equation*}
    \label{Lipschitz R}
\end{assumption}

\begin{assumption}
    The function $L(x,\gamma): \mathcal{X} \times \mathcal{Y} \to \mathbb{R}$ is $\rho_L$-smooth. That is, $\exists \; \rho_L >0$, $\forall \; (x_1,\gamma_1) \in \mathcal{X} \times \mathcal{Y}$ and $(x_2,\gamma_2) \in \mathcal{X} \times \mathcal{Y}$ we have
    \begin{equation*}
        \norm{\nabla L(x_1,\gamma_1)-\nabla L(x_2,\gamma_2)} \leq \rho_L \norm{(x_1,\gamma_1) -(x_2,\gamma_2)}.
    \end{equation*}
    \label{hdh}
\end{assumption}

With these assumptions in place, we can now present the main theorem of this subsection.

\begin{theorem}
 Under Assumptions \ref{saddle assumption}-\ref{subg} and \ref{steps}-\ref{hdh},  for any $\epsilon > 0$, there exists $\mathbf{B} > 0$ such that if $\norm{b_x(n)}$ and $\norm{b_\gamma(n)} \leq \mathbf{B}$, then  $ \exists \; n_0 \in \mathbb{N}$ such that  $ \forall \; n \geq n_0$, we have $(x(n),\gamma(n)) \in N_\epsilon (\mathcal{S})$ a.s.
\end{theorem}
\begin{proof}
As in Theorem \ref{almstsure} , we focus on the trajectory $(\Hat{x}(t),\Hat{\gamma}(t))$ instead of $(x(n),\gamma(n))$ for precisely the same reason.

Let us first define $V(x,\gamma) : \mathcal{X} \times \mathcal{Y} \to \mathbb{R}^{+}$ for each $(x^\ast,\gamma^\ast) \in \mathcal{S}$  as
\begin{equation*}
    V(x,\gamma) = \mathbb{D}_{R_x} (x^\ast,x) + \mathbb{D}_{R_\gamma} (\gamma^\ast,\gamma)
\end{equation*}

and $V^\ast(x,\gamma) : \mathcal{X} \times \mathcal{Y} \to \mathbb{R}^{+}$ as
\begin{equation*}
    V^\ast(x,\gamma) = \mathbb{D}_{R_x} (x^\ast,x) + \mathbb{D}_{R_\gamma} (\gamma^\ast,\gamma).
\end{equation*}
It is already shown in the proof of Theorem \ref{almstsure} that $V^\ast(x,\gamma)$ is continuous in $\mathcal{X} \times \mathcal{Y}$. Furthermore, since $\mathcal{X} \times \mathcal{Y}$ is assumed to be compact, $V^\ast (x,\gamma)$ is uniformly continuous. This implies that $\forall \; \epsilon > 0$, $\exists \; \delta_1 > 0$ such that if $\norm{(x_1,\gamma_1)-(x_2,\gamma_2)} \leq \delta_1$, $\forall \; (x_1,\gamma_1), (x_2,\gamma_2) \in \mathcal{X}\times \mathcal{Y}$, then $\norm{V^\ast(x_1,\gamma_1)-V^\ast(x_2,\gamma_2)} \leq \frac{\epsilon}{6}$.

Consider $(x(t),\gamma(t))$ as any Caratheodory solution of the projected dynamical system \eqref{projected primalip1}-\eqref{projected dualip1} with $b_x(t),b_\gamma(t) = 0$ and any initial condition $(x_0,\gamma_0) \in \mathcal{X} \times \mathcal{Y}$. Since $\mathcal{X} \times \mathcal{Y}$ is compact, as proven in step $1$ of Theorem \ref{almstsure}, there exists $T > 0$ such that $V(x(T),\gamma(T)) \leq \frac{\epsilon}{6}$ regardless of the initial condition. 
%Refer this $T$ for the remainder of the proof.

 From the proof of Theorem \ref{almstsure}, it is clear that $\forall \; \delta > 0$, $\exists \; t_0(T,\delta)$ such that $\forall \; t \geq t_0$ we have  (a.s.)
   \begin{equation}
       V(\Hat{x}(t+T),\Hat{\gamma}(t+T)) \leq V(x_t(T),\gamma_t(T)) + \rho_R \delta^2 + \rho_R \delta (\norm{x_t(T)-x^\ast} + \norm{\gamma_t(T) -\gamma^\ast}) \leq V(x_t(T),\gamma_t(T)) + \rho_R \delta^2 + \rho_R\delta D.
       \label{720}
   \end{equation}
  Here, $(x_t(s),\gamma_t(s))$ represents the Caratheodory solution of the projected dynamical system \eqref{projected primalip1}-\eqref{projected dualip1} with the initial condition $(x_0,\gamma_0) = (\Hat{x}(t),\Hat{\gamma}(t))$ and $D$ is the diameter of the set $\mathcal{X}$ and $\mathcal{Y}$.
Let $(x_t^1(s),\gamma_t^1(s))$ denote the Caratheodory solutions of the projected dynamical system \eqref{projected primalip1}-\eqref{projected dualip1} with the same initial condition, but with $b_x(t),b_\gamma(t) = 0$. Choose $\delta$ in such a way that $\rho_R \delta^2 + \rho_R\delta D \leq \frac{2 \epsilon}{3}$.

   Consider the following term 
   \begin{equation*}
       \begin{split}
           & \frac{d}{ds}(\norm{x_t(s)-x_t^1(s)}_{x_t^1(s)}^2 + \norm{\gamma_t(s)-\gamma_t^1(s)}_{\gamma_t^1(s)}^2 ).
       \end{split}
   \end{equation*}
   
   Next, we focus on the term involving only $x_t$ and $x_t^1$
   \begin{equation}
       \begin{split}
           & \frac{d}{ds}(\norm{x_t(s)-x_t^1(s)}_{x_t^1(s)}^2)
           \\ = & \frac{d}{ds}( (x_t(s)-x_t^1(s))^\top \nabla^2 R_x(x_t^1(s))(x_t(s)-x_t^1(s)))
           \\ = & ( (x_t(s)-x_t^1(s))^\top \nabla^2 R_x(x_t^1(s))(\Dot{x}_t(s)-\Dot{x}_t^1(s))) + ( (x_t(s)-x_t^1(s))^\top  \frac{d}{ds}(\nabla^2 R_x(x_t^1(s)))(x_t(s)-x_t^1(s))).
       \end{split}
       \label{1000}
   \end{equation}
  
   Consider the second term on the right-hand side of equation \eqref{1000}
   \begin{equation}
       ( (x_t(s)-x_t^1(s))^\top  \frac{d}{ds}(\nabla^2 R_x(x_t^1(s)))(x_t(s)-x_t^1(s))) \leq \norm{\frac{d}{ds}(\nabla^2 R_x(x_t^1(s)))} \norm{x_t(s)-x_t^1(s)}^2.
       \label{dg}
   \end{equation}
  
   Hence,
   \begin{equation}
       \begin{split}
           & \norm{ \frac{d}{ds}(\nabla^2 R_x(x_t^1(s)))}
           \\ = & \norm{\lim\limits_{h \to 0^+} \frac{\nabla^2 R_x (x_t^1(s+h))-\nabla^2 R_x(x_t^1(s))}{h}}
           \\ = &\lim\limits_{h \to 0^+} \frac{\norm{\nabla^2 R_x(x_t^1(s+h))-\nabla^2 R_x(x_t^1(s))}}{h}
           \\ \leq & K_2 \lim\limits_{h \to 0^+} \frac{\norm{x_t^1(s+h)-x_t^1(s)}}{h}
           \\ = & K_2 \norm{\lim\limits_{h \to 0^+}\frac{x_t^1(s+h)-x_t^1(s)}{h}} 
           = K_2 \norm{\frac{d}{ds} x_t^1(s)} \leq K_2 K_3 \norm{\nabla_x L(x_t^1(s),\gamma_t^1(s))} \leq K_2K_3 G.
       \end{split}
       \label{74}
   \end{equation}
   The first inequality in \eqref{74} arises in view of Assumption \ref{Lipschitz R}. The second inequality stems from the fact that $(x_t^1(s),\gamma_t^1(s))$ is a Caratheodory solution of the ordinary differential inclusions given in \eqref{projected primalip1}-\eqref{projected dualip1} with $b_x(t),b_\gamma(t) = 0$. However, in view of lemma \ref{equivalence}, $(x_t^1(s),\gamma_t^1(s))$ is also a Caratheodory solution of the following projected dynamical system with initial condition $x_t^1(0) = \Hat{x}(t)$ and $\gamma_t^1(0) = \Hat{\gamma}(t)$. 
   \begin{equation*}
       \Dot{x}(s) = \mathcal{P}_{\mathcal{T}_\mathcal{X}(x)}^{x(s)} (-\nabla_x L(x(s), \gamma(s))) \; \; \text{and} \; \; \Dot{\gamma}(s) = \mathcal{P}_{\mathcal{T}_\mathcal{Y}(\gamma)}^{\gamma(s)} (\nabla_\gamma L(x(s), \gamma(s))).
   \end{equation*}
    Noting that $0 \in \mathcal{T}_\mathcal{X}(x(t))$), we have
\begin{equation*}
\norm{\Dot{x}_t^1(s) + \nabla_x L(x_t^1(s),\gamma_t^1(s)) }_{x_t^1(s)} \leq \norm{\nabla_x L(x_t^1(s),\gamma_t^1(s))}_{x_t^1(s)}.
\end{equation*}
Hence,
\begin{equation*}
    \norm{\Dot{x}_t^1(s) }_{x_t^1(s)} \leq 2 \norm{\nabla_x L(x_t^1(s),\gamma_t^1(s))}_{x_t^1(s)}. 
\end{equation*}
In view of norm equivalence in finite dimension, we get $\exists \; K_3(x_t^1(s)) > 0 $
\begin{equation*}
    \norm{\Dot{x}_t^1(s) } \leq K_3 (x_t^1(s))\norm{\nabla_x L(x_t^1(s),\gamma_t^1(s))} \leq K_3 (x_t^1(s)) G . 
\end{equation*}
 However, $K_3$ is independent of $x_t^1(s)$ in view of lemma \ref{normequivalence}.

Hence, from \eqref{dg}, we obtain
\begin{equation}
     ( (x_t(s)-x_t^1(s))^\top  \frac{d}{ds}(\nabla^2 R_x(x_t^1(s)))(x_t(s)-x_t^1(s))) \leq K_2K_3G \norm{x_t(s)-x_t^1(s)}^2. 
     \label{bomb}
\end{equation}
Now, let's consider the first term of the RHS of \eqref{1000}: 
\begin{equation}
    \begin{split}
        & ( (x_t(s)-x_t^1(s))^\top \nabla^2 R_x(x_t^1(s))(\Dot{x}_t(s)-\Dot{x}_t^1(s))) 
        \\ = & (x_t^1(s)-x_t(s))^\top \nabla^2 R(x_t^1(s)) (\nabla^2 R_x(x_t^1(s))^{-1} (-\nabla_x L(x_t^1(s),\gamma_t^1(s)) -\eta_1^1(s) +\nabla_x L(x_t(s),\gamma_t(s)) + \eta_1(s)-b_x(s))  
        \\  & +  (x_t^1(s)-x_t(s))^\top \nabla^2 R(x_t^1(s)) \{ (\nabla^2 R_x(x_t^1(s))^{-1} (-\nabla_x L(x_t(s),\gamma_t(s)) - \eta_1(s)+b_x(s)) - \Dot{x}_t(s) \}
    \end{split}
    \label{mnss}
\end{equation}
Equation \eqref{mnss} holds since $(x_t(s),\gamma(s))$ is a Caratheodory solution of the ordinary differential inclusions  given in \eqref{projected primalip1}-\eqref{projected dualip1} with initial condition $(x_t(0),\gamma_t(0)) = (\Hat{x}(t),\Hat{\gamma}(t))$. Consequently,  there exists $\eta_1(s) \in \mathcal{N}_\mathcal{X}(x_t(s))$ such that $\Dot{x}_t(s) = \nabla^2 R_x(x_t(s))^{-1} (-\nabla_x L(x_t(s),\gamma_t(s)) - \eta_1(s)+b_x(s))$ almost everywhere. 

With a similar line of  reasoning, $\exists \; \eta_1^1(s) \in \mathcal{N}_\mathcal{X}(x_t^1(s))$ such that $\Dot{x}_t^1(s) =\nabla^2 R_x(x_t^1(s))^{-1}(- \nabla_x L(x_t^1(s),\gamma_t^1(s)) -\eta_1^1(s)) $ a.e.

Considering only the second term on the RHS of \eqref{mnss}, we have
\begin{equation}
    \begin{split}
        & (x_t^1(s)-x_t(s))^\top \nabla^2 R(x_t^1(s)) \{ (\nabla^2 R_x(x_t^1(s))^{-1} (-\nabla_x L(x_t(s),\gamma_t(s)) - \eta_1(s)+b_x(s)) - \Dot{x}_t(s) \}
        \\ = &  (x_t^1(s)-x_t(s))^\top (-\nabla_x L(x_t(s),\gamma_t(s)) -\eta_1(s) + b_x(s) - \nabla^2 R_x(x_t^1(s))\Dot{x}_t(s) )
        \\ = &  (x_t^1(s)-x_t(s))^\top (\nabla^2 R_x(x_t(s)) \Dot{x}_t(s) - \nabla^2 R_x(x_t^1(s)) \Dot{x}_t(s))
        \\ \leq & \norm{x_t^1(s)-x_t(s)} \norm{\nabla^2 R_x (x_t(s)) - \nabla^2 R_x (x_t^1(s))} \norm{\Dot{x}_t(s)}
        \\ \leq & K_2 \norm{x_t^1(s)-x_t(s)}^2 \norm{\Dot{x}_t(s)} \leq K_2 K_3 \norm{x_t^1(s)-x_t(s)}^2 (\norm{\nabla_x L(x_t(s),\gamma_t(s))}+ \norm{b_x(s)}) \\ \leq &   K_2 K_3 \norm{x_t^1(s)-x_t(s)}^2 (G + \mathrm{B})
    \end{split}
    \label{jumd}
\end{equation}
The second inequality in \eqref{jumd} arises from the Lipschitz continuity of the function $\nabla^2R_x(x)$ and using the fact $\norm{b_x(S)} \leq \mathrm{B}$. 
The fourth inequality stems from the fact that $(x_t(s),\gamma_t(s))$  is a Caratheodory solution of \eqref{projected primalip1}-\eqref{projected dualip1}. 
%as explained earlier in this proof.

Now, let's examine the first term on the RHS of \eqref{mnss}
\begin{equation*}
    \begin{split}
        & (x_t^1(s)-x_t(s))^\top  (-\nabla_x L(x_t^1(s),\gamma_t^1(s)) -\eta_1^1(s) +\nabla_x L(x_t(s),\gamma_t(s)) + \eta_1(s)+b_x(s)) 
        \\ \leq & (x_t^1(s)-x_t(s))^\top(-\nabla_x L(x_t^1(s),\gamma_t^1(s))+\nabla_x L(x_t(s),\gamma_t(s))) + (\norm{x_t^1(s)-x_t(s)}^2+1)\mathrm{B}. 
    \end{split}
\end{equation*}
Note that $(x_t^1(s)-x_t(s))^\top(\eta_1(s)-\eta_1^1(s)) \leq 0 $ due the definition of normal cone. 

Also, $(x_t^1(s)-x_t(s))^\top b_x(s) \leq (\norm{x_t(s)-x_t^1(s)}^2 +1) \mathrm{B}$. Assume that $\norm{b_x(s)},\norm{b_\gamma(s)} \leq \mathrm{B}$ $\forall \; s >0$. 

Hence, adding equations \eqref{mnss} and \eqref{bomb}, we obtain
\begin{equation}
    \begin{split}
         \frac{d}{ds} (\norm{x_t(s)-x_t^1(s)}_{x_t^1(s)}^2) \leq & (K_2K_3G  + (G+\mathrm{B}) K_2K_3 + \mathrm{B}) \norm{x_t^1(s)-x_t(s)}^2 + \mathrm{B} \\ & + (x_t^1(s)-x_t(s))^\top (\nabla_x L(x_t(s),\gamma_t(s))-\nabla_x L(x_t^1(s),\gamma_t^1(s)))
    \end{split}
    \label{x variable}
\end{equation}
With the similar procedure, we can show that 
\begin{equation}
    \begin{split}
         \frac{d}{ds} (\norm{\gamma_t(s)-\gamma_t^1(s)}_{\gamma_t^1(s)}^2)  & \leq  (K_2K_3G  + (G+\mathrm{B}) K_2K_3 + \mathrm{B}) \norm{\gamma_t^1(s)-\gamma_t(s)}^2 + \mathrm{B} \\ & + (\gamma_t^1(s)-\gamma_t(s))^\top (\nabla_\gamma L(x_t(s),\gamma_t(s))-\nabla_\gamma L(x_t^1(s),\gamma_t^1(s))).
    \end{split}
    \label{y variable}
\end{equation}
Notice that due to the $\rho_L$-smoothness of the function $L(x,\gamma)$, we have
\begin{equation*}
    \begin{split}
        & (x_t^1(s)-x_t(s))^\top (\nabla_x L(x_t(s),\gamma_t(s))-\nabla_x L(x_t^1(s),\gamma_t^1(s))) + (\gamma_t^1(s)-\gamma_t(s))^\top (\nabla_\gamma L(x_t(s),\gamma_t(s))-\nabla_\gamma L(x_t^1(s),\gamma_t^1(s)))
        \\ \leq & \rho_L(\norm{x_t^1(s)-x_t(s)}^2 + \norm{\gamma_t(s)-\gamma_t^1(s)}^2).
    \end{split}
\end{equation*}
Denoting $a = K_2K_3G + (G+\mathrm{B}) + \mathrm{B}+ \rho_L$, adding both \eqref{x variable} and \eqref{y variable} and in view of lemma \ref{normequivalence} we have, $\exists \; K_4 > 0$ such that 
\begin{equation*}
    \begin{split}
         \frac{d}{ds} (\norm{x_t(s)-x_t^1(s)}^2+\norm{\gamma_t(s)-\gamma_t^1(s)}^2)  \leq a K_4  (\norm{x_t(s)-x_t^1(s)}^2+\norm{\gamma_t(s)-\gamma_t^1(s)}^2) + \mathrm{B} K_4 .
    \end{split}
\end{equation*}
Since $x_t(0)=x_t^1(0)$ and $\gamma_t(0) = \gamma_t^1(0)$, for all $0 \leq s \leq T$, we have 
\begin{equation}
    \norm{x_t(s)-x_t^1(s)}^2+\norm{\gamma_t(s)-\gamma_t^1(s)}^2 \leq aK_4 \int\limits_{0}^{s} \{ \norm{x_t(z)-x_t^1(z)}+\norm{\gamma_t(z)-\gamma_t^1(z)}^2 \} dz  + \mathrm{B}K_4 T.
\end{equation}
Applying Gronwall inequality yields
\begin{equation}
     \norm{x_t(s)-x_t^1(s)}^2+\norm{\gamma_t(s)-\gamma_t^1(s)}^2 \leq \mathrm{B} K_4 T e^{a K_4T}.
     \label{rumpaii}
\end{equation}
From \eqref{rumpaii}, it follows that $\exists \; \mathbf{B} > 0$ such that if $\mathrm{B} \in (0, \mathbf{B})$, then
\begin{equation*}
    \norm{(x_t(T),\gamma_t(T))-(x_t^1(T),\gamma_t^1(T))} \leq \delta_1
\end{equation*}
This implies 
\begin{equation*}
    \norm{V^\ast((x_t(T),\gamma_t(T))-V^\ast(x_t^1(T),\gamma_t^1(T))} \leq \frac{\epsilon}{6}.
\end{equation*}
However, since $V^\ast(x_t^1(T),\gamma_t^1(T)) \leq \frac{\epsilon}{6}$ due to $(x_t^1(0),\gamma_t^1(0)) \in \mathcal{X} \times \mathcal{Y}$, it follows that $V^\ast(x_t(T),\gamma_t(T)) \leq \frac{\epsilon}{3}$. Thus, from \eqref{720}, $\exists$ $t_0 > 0$ such that $\forall$ $t \geq t_0$, $\exists \; (x^\ast,\gamma^\ast) \in \mathcal{S}$ such that $V(\Hat{x}(t+T),\Hat{\gamma}(t+T)) \leq \epsilon$ a.s.

From this point, the conclusion of the Theorem can be established by utilizing the  inequality: 
$V(x,
\gamma) \geq \frac{\rho_R}{2}(\norm{x-x^\ast}^2 + \norm{\gamma-\gamma^\ast}^2)$. This inequality is valid due to the $\rho_R$-strong convexity of the functions $R_x$ and $R_\gamma$.
\end{proof}
The following Corollary immediately follows from the above Theorem. 
\begin{Corollary}
    When $\mu(n) \to 0$ in \eqref{approx}, the iterates of the   SZSPMD algorithm converge almost surely to the set $\mathcal{S}$.
\end{Corollary}

However, if we reduce 
$\mu(n)$ to very small values, the impact of system noise could dominate the difference between function values in \eqref{approx}. In such cases, \eqref{approx} may fail to accurately represent the approximated gradient of the function 
$L(x,\gamma)$ (refer to \cite{9186148} for more details). Thus, it's crucial to carefully select $\mu(n)$. However, the specific methodology for choosing 
$\mu(n)$ is beyond the scope of this paper. For further insights on selecting 
$\mu(n)$, we  refer the reader to \cite{10143924}, although it approaches the topic from a different perspective, focusing on concentration bounds rather than the dynamic viewpoint. 
\section{Conclusion and Future Direction}
The work in this paper has presented novel methodologies for addressing saddle point problem, leveraging projected dynamical systems in non-Euclidean spaces. 
The proposed approach exploits the geometric characteristics of saddle point problems, thereby expanding the scope to address a wider range of complex problems beyond traditional Euclidean frameworks.

We have demonstrated the theoretical foundations of these dynamical systems, establishing the existence of viable  Caratheodory solutions and analyzing their stability properties. 
Furthermore, we discretize the proposed dynamical system and demonstrate its asymptotic equivalence with the saddle point mirror descent algorithm. Both stochastic zeroth and first-order saddle point mirror descent algorithms are discussed, and the almost sure convergence analysis of their iterates is presented, leveraging the stability analysis of the projected dynamical system in a non-Euclidean domain.   In summary, we have laid the theoretical groundwork for various practical implementations of projected dynamical systems in non-Euclidean domains across diverse scenarios.
 However, it's important to note that the analysis in this paper is confined to the asymptotic analysis of SSPMD algorithms and SZSPMD   algorithms.
 Looking forward, future research will delve into the finite time analysis of these algorithms through determining a  
  concentration inequalities for the error between primal and dual updates.
 Additionally, an important extension would be to extend the projected dynamical system within a non-Euclidean domain to tackle non-convex-concave saddle point problems. By continuing to explore these avenues, we can further enhance the applicability and effectiveness of projected dynamical systems in Non-Euclidean domain.
\bibliographystyle{IEEEtran}
\bibliography{ref}
\end{document}